\documentclass[final,3p,times,12pt]{elsarticle}

\usepackage{amssymb}
\usepackage{amsthm}
\usepackage{amsmath}
\usepackage{latexsym}
\usepackage{times}
\usepackage{graphicx}
\usepackage{epsfig}
\usepackage{epstopdf}
\usepackage{mathrsfs}
\usepackage{bm} 
\usepackage{appendix} 

\usepackage[active]{srcltx}
\usepackage{subcaption}

\usepackage{adjustbox}

\usepackage{hyperref}  		
\hypersetup{colorlinks=true,citecolor={red},linkcolor={blue},linktocpage} 
\usepackage{color}
\definecolor{armygreen}{rgb}{0.29, 0.33, 0.13}
\definecolor{maroon(x11)}{rgb}{0.69, 0.19, 0.38}
\definecolor{darkorange}{rgb}{1.0, 0.55, 0.0}
\definecolor{auburn}{rgb}{0.43, 0.21, 0.1}

\newcommand{\uinc}{u_{\rm inc}}

\newcommand{\nl}{n_{\lambda}}

\newcommand{\Sph}{\mathbb{S}}

\newcommand{\vect}[1]{\boldsymbol{\mathbf{#1}}} 

\usepackage{color}
\definecolor{hotcolor}{rgb}{1,0,0}

\journal{Elsevier}

\begin{document}

\begin{frontmatter}




\title{\textbf{Highly accurate acoustic scattering: Isogeometric Analysis
coupled with local high order Farfield Expansion ABC}}


\author[label1]{Tahsin Khajah}
\ead{tkhajah@uttyler.edu}

\author[label2]{Vianey Villamizar\corref{cor1}}
\ead{vianey@mathematics.byu.edu}

\address[label1]{Department of Mechanical Engineering, University of Texas at Tyler, Texas}
\address[label2]{Department of Mathematics, Brigham Young University, Provo, UT}
\cortext[cor1]{Corresponding author}

\begin{abstract}
This work is concerned with a unique combination of high order local absorbing boundary conditions (ABC) with a general curvilinear Finite Element Method (FEM) and its implementation in Isogeometric Analysis (IGA) for time-harmonic acoustic waves. The ABC employed were recently devised 
by Villamizar, Acosta and Dastrup  [J. Comput. Phys. 333  (2017) 331] .  They are derived from exact Farfield Expansions representations of the outgoing waves in the exterior of the regions enclosed by the artificial boundary.  As a consequence, the error due to the ABC on the artificial boundary can be reduced conveniently such that the dominant error comes from the volume discretization method used in the interior of the computational domain. Reciprocally, the error in the interior can be made as small as the error at the artificial boundary by appropriate implementation of {\it p-} and {\it h}- refinement.
We apply this novel method to cylindrical, spherical and arbitrary shape scatterers including a prototype submarine.
Our numerical results exhibits spectral-like approximation and high order convergence rate.
Additionally, they show that the proposed method can reduce both the pollution and artificial boundary errors to negligible levels even in very low- and high- frequency regimes with rather coarse discretization densities in the IGA. As a result, we have developed a highly accurate computational platform to numerically solve time-harmonic acoustic wave scattering in two- and three-dimensions. 

\end{abstract}


\end{frontmatter}


\section{Introduction} \label{Section.Intro}

The development of efficient, robust, and easy to implement numerical methods for exterior acoustic scattering problems has been intensively studied. In spite of these efforts, challenges still remain for the different approaches followed. For example finite element methods (FEM), which have become very popular for its ability to model complex geometries and its mathematical robustness, suffer from dispersion errors \cite{Babuska-Sauter1997,Ihlenburg}. Therefore,  large amount of computational resources are required to obtain accurate approximations. Another challenge of 
traditional FEM is to reduce the induced errors introduced by the approximated representation of the boundary geometry \cite{Hughes-Cottrell2005,Nguyen-Bordas}. Additionally, all volume discretization methods applied to exterior acoustic scattering require to introduce an artificial boundary to truncate the unbounded physical domain and to impose 
an absorbing boundary condition (ABC) on it \cite{GivoliReview2,Tsynkov1998}. As a consequence, the approximate solution is affected to some degree by spurious reflection  from the boundary.  

Among recent efforts made to overcome the above challenges, we find the work by Turkel et al. in \cite{Turkel-Farhat2004}. They reduced the Helmholtz equation to a new one without the main oscillatory term. This equation was combined with the well-known ABC named  Bayliss-Gunzburger-Turkel: BGT-1 and BGT-2 \cite{Bayliss01}. In two dimensions, they also considered another ABC defined from the second order operator (BGTH) that annihilates the leading order term in Karp's expansion \cite{Karp}, as shown in \cite{Grote-Keller01}. Then, they applied linear finite elements to approximate  scattering problems. They found that the non-oscillatory equation produces slightly better results than the Helmholtz equation when both are combined with BGT-2 for the low frequencies $k=3,5$. However for very low frequencies such as $k=0.01$, the Helmholtz equation combined with BGTH outperforms the non-oscillatory equation coupled with any of their three ABCs by several orders of magnitude. This was previously observed by Grote and Keller in  \cite{Grote-Keller01}. However, in both work their results were limited by the low order of the finite element basis (linear) and the low order of the ABC employed. 

In \cite{Antoine2009}, Kechroud et al. considered two-dimensional acoustic scattering from circular, elliptical and a submarine-like shaped scatterers in two dimensions. 
An enriched FEM with plane waves \cite{Laghrouche2000,Laghrouche2002} was employed  in an effort to avoid dispersion. This method was coupled with a Pad\'e-type ABC. They conducted experiments for a circular scatterer with frequencies $k=1,5,10,20$ and compared against the exact solution for the bounded problem with the Pad\'e-type ABC. 
They obtained accurate results in the order of $10^{-4}$ for the relative $L^2$ error both on the computational domain and at the artificial boundary. A virtue of this approach is its ability to reduce the required number of elements per wavelength which in turn reduces the degrees of freedom used. By doing this, they found good approximations even for moderately high frequencies such as $k=60$. Unfortunately, this technique is limited to two-dimensional problems and suffers from the inaccurate approximation of the scatterer boundaries.

In an effort to overcome the errors associated to a poor representation of the boundary geometry, Khajah et al. \cite{Tahsin2016} applied Isogeometric Finite Element Analysis (IGA) to the exterior acoustic scattering from a circular scatterer in two dimensions using BGT-2 as absorbing boundary condition. They compared their numerical solutions against the analytical solution obtained for the bounded problem with the BGT-2 absorbing boundary condition. As a result, they eliminated the errors due to the ABC. They reported errors in the order of $10^{-3}$ for high frequencies such as $k=200, 500$. These authors also employed a Pad\'e type ABC in \cite{Tahsin_JCA2018}  that slightly improved the  artificial boundary errors incurred by using  BGT2.

Similarly, Dinachandra and Raju \cite{Dinachandra-Raju2018} implemented an IGA technique with plane wave enrichment that they called Partition of Unity Isogeometric Analysis (PUIGA). They applied it to several benchmark problems including the exterior acoustic scattering of a plane wave from a circular cylinder. They considered two boundary value problems (BVPs). One of them   
contained a non-homogeneous Sommerfeld boundary condition that was obtained from substituting the sound-hard exact solution for the scattering problem into the Sommerfeld condition. The other one employed the BGT-2 absorbing boundary condition. Their numerical solutions for both problems were compared against the sound-hard scattering exact solution. 
For the first problem without the domain truncation error, they obtained approximations in the order of $10^{-4}$ for $k=10$ with $q=8$ plane waves; and $10^{-5}$ for $k=10$ with $q=16$. The merit of this technique is that these errors were attained for relatively coarse grids. For the second problem with the BGT-2 condition, the errors obtained with this technique were two order higher than those of the first problem. Clearly, the use of a low order BGT-2 absorbing condition caused the deterioration of the numerical solution. 

Schmidt and Heier \cite{Schmidt2015} also found numerical solutions for the plane wave scattering from a circular cylinder. These authors used high order FEM coupled with  Feng's absorbing Boundary conditions of several orders at the artificial boundary. They studied the convergence to the exact solution 
with respect to the radius of the artificial boundary. The authors were able to obtain very sharp theoretical estimates for the $L^2$ error. These estimates were numerically verified for large values of the radius $R$ of the artificial boundary.  These computations were only carried out for a fixed frequency $k=1$.

There has been other attempts using high order method for the approximation of the Helmholtz equation such as plane wave enriched finite element and IGA techniques. The reader is referred to the article \cite{Dinachandra-Raju2018} for a rather complete set of these contributions . Unfortunately to the best of our knowledge, most of this work has been done in two dimensions and without employing high order local ABC. As a consequence, the high order approximation of the interior methods is negatively affected by the low order approximation of the ABC which results in an overall low order method.

In this work, we propose a numerical method consisting of an IGA technique employing arbitrary high order NURBS bases combined with recently developed high order local Farfield Expansions absorbing boundary conditions \cite{JCP2017}. 
These conditions are defined  from truncated versions of exact series representations of the outgoing waves outside a circular and spherical artificial boundary in 2D and 3D, respectively. The series representations employed are Karp's farfield expansion \cite{Karp} in 2D, and Wilcox's farfield expansion \cite{Wil-1956} in 3D. Therefore, the resulting absorbing boundary conditions called
Karp's farfield expansion ABC (KFE) and Wilcox farfield expansion ABC (WFE), respectively, can be considered exact absorbing boundary conditions.  The angular functions appearing in Wilcox's or Karp's farfield expansions are part of the unknowns. To determine these angular functions, the recurrence formulas derived from Wilcox's or Karp's theorems which do not disturb the local character of the ABC are employed. Moreover, the order of the error at the boundary induced by this ABC can be easily reduced by simply adding as many terms as needed to the truncated farfield expansions.

\section{Formulation of the problem}
\label{Section.Formulation}


We consider the scattering of a time-harmonic incident wave from a single obstacle in two or three dimensions. 
The incident wave is a monochromatic plane wave, $u_{inc}(\textbf{x})e^{-i\omega
t}= e^{i k \textbf{x} \cdot \textbf{d}}e^{-i\omega t}$, where
$\textbf{d}$ is a unit vector that points in the \emph {direction of
incidence,} and $i=\sqrt{-1}$. This incident wave is impinging upon
 an impenetrable obstacle that occupies a simply connected bounded region  with boundary $\Gamma$. The open unbounded region in the exterior of $\Gamma$ is denoted as 
$\Omega^*$. The incident field $\uinc$  satisfies the Helmholtz equation in $\Omega^*$, and the scattered field 
$u$ solves the boundary value problem (BVP):
\begin{eqnarray}
&& \Delta u + k^2 u = f \quad\qquad \text{in $\Omega^*$}, \label{BVPsc1} \\
&& Z \partial_n u + (1-Z) u = -\Big(Z
\partial_n u_{inc}+ (1-Z) u_{inc}\Big) \qquad\qquad \text{on $\Gamma$,}  \label{BVPsc2} \\
&& \lim_{r \rightarrow \infty} r^{(\delta-1)/2} \left( \partial_{r} u- \mathrm{i} k u \right) = 0.\label{BVPsc3}
\end{eqnarray}
 The wave number $k$ and the source $f$ may vary in space. Equation (\ref{BVPsc3}) is known as the Sommerfeld radiation condition where $r = |\textbf{x}|$ and $\delta=2$ or 3 for two or three dimensions, respectively. It implies that $u$ is an outgoing wave. In equation (\ref{BVPsc2}), $Z=0$ or $Z=1$.
If $Z=0$ the boundary condition (\ref{BVPsc2} ) is a Dirichlet condition that models
acoustically soft obstacles and if $Z=1$ is a Neumann condition which models hard obstacles. We  consider both BVPs in this work.
 These boundary value problems are well-posed under classical and weak formulations \cite{ColtonKress02,Nedelec01,McLean2000}. 
 The analysis and numerical computation in this article can be easily extended to Robin boundary conditions, and to a bounded penetrable scatterer with inhomogeneous and anisotropic properties. 

As pointed out in the introduction, the unbounded BVP (\ref{BVPsc1})-(\ref{BVPsc3}) needs to be transformed into a bounded BVP before a numerical solution can be sought. In Villamizar et al. \cite{JCP2017}, this transformation was carried out by introducing
a circular (two dimensions) and a spherical (three dimensions) artificial boundaries, respectively. This was followed by defining high order local absorbing boundary conditions based on farfield expansions ABC on these artificial boundaries.


\section{The scattering problem in two dimensions. Karp's expansion}
\label{2DProblem}

For the two-dimensional case, an equivalent problem to (\ref{BVPsc1})-(\ref{BVPsc3}) was introduced in \cite{JCP2017} whose truncated version is given by
 \begin{align}
& \Delta u + k^2 u = f, \quad\qquad  && \text{in $\Omega$}, \label{BVPBd1} \\
&Z {\partial_n u} + (1-Z) u = -\Big(Z
\partial_n u_{inc}+ (1-Z) u_{inc}\Big)&& \text{on $\Gamma$,}  \label{BVPBd2} \\
& u(R,\theta)=H_0(kR) \sum_{l=0}^{L-1} \frac{F_l(\theta)}{(kR)^l} + H_1(kR)\sum_{l=0}^{L-1} \frac{G_l(\theta)}{(kR)^l},\label{BVPBd3} \\
&\partial_{r} u(R,\theta) = \partial_{r}\left( H_0(kr) \sum_{l=0}^{L-1} \frac{F_l(\theta)}{(kr)^l} + H_1(kr)\sum_{l=0}^{L-1} \frac{G_l(\theta)}{(kr)^l}\right) \bigg|_{r=R},
\label{BVPBd4} \\
& \partial_{r}^2 u(R,\theta) = \partial_{r}^2\left( H_0(kr)\sum_{l=0}^{L-1} \frac{F_l(\theta)}{(kr)^l} + H_1(kr)\sum_{l=0}^{L-1} \frac{G_l(\theta)}{(kr)^l}\right) \bigg|_{r=R}, \label{BVPBd5}\\
& 2 l F_{l}(\theta) = - l^2 G_{l-1}(\theta) - d^2_{\theta} G_{l-1}(\theta), \qquad && \text{for $l=1,2, \dots$}. \label{Recurrence2}\\
& 2 l G_{l}(\theta) = (l-1)^2 F_{l-1}(\theta) + d^2_{\theta}  F_{l-1}(\theta) , \qquad && \text{for $l=1,2, \dots$} \label{Recurrence1}
\end{align}
where $R$ is the radius of a circular artificial boundary $S_R$ enclosing the scatterer and $\Omega$ is the annular region bounded by $\Gamma$ and $S_R$. 
The equations (\ref{BVPBd3})-(\ref{BVPBd5}) for the \textit{double} family of unknown farfield functions $F_l$ and $G_l$, supplemented by the recurrence formulas (\ref{Recurrence2})-(\ref{Recurrence1}), constitute the novel \textit{Karp's Farfield Expansion}  absorbing boundary condition $(\text{KFE})$ that was proposed in \cite{JCP2017}. Notice, that the set of equations (\ref{BVPBd3})-(\ref{BVPBd5}) are enough to determine the approximations of  $u$, $F_0$ and $G_0$ at the artificial boundary. The last two recurrence formulas serve to determine the angular functions $F_l$ and $G_l$, for $l=1\dots L-1$.
In what follows, we will assume that $f=0$ for simplicity.
It was shown in \cite{JCP2017} that the numerical solution of  (\ref{BVPBd1})-(\ref{Recurrence1}) exhibits second order convergence to the exact solution if a standard second order finite difference method is employed in the interior of the computational domain. These results were obtained even for cases where the artificial boundary was imposed extremely close to the scatterer (see Fig. 2 in \cite{JCP2017}). The number of terms employed by the KFE was relatively small
(usually three to eight) in many practical situations.

One of the main purpose in this article is to further exploit the high order property of the  KFE by coupling them with a high order isogeometric finite element method. We will show that it is possible to obtain orders of convergence grater than second order and high accuracy in the numerical solutions by  appropriately adjusting the farfield expansion number of terms and by performing {\it h-} and {\it p}- refinements.

\subsection{Weak formulation and finite element approximation.}
\label{Weakformulation}

We will derive the weak form of the BVP  (\ref{BVPBd1} )-(\ref{Recurrence2}) for the Dirichlet case ($Z=0$) by defining the function spaces:
\begin{align}
&\mathscr{S}=\{(u,F_0,G_0\dots F_{L-1},G_{L-1} ) \,|\,u=-u_{inc}  \mbox{ on } \Gamma,\, u\in H^1(\Omega),\, 
F_l,G_l\in H^1(S_R),\, l=0\dots L-2, \nonumber\\
&\qquad\,\,\,F_{L-1}, G_{L-1} \in H^0(S_R)\}\\
&\mathscr{S}_0=\{v \,|\,v=0  \mbox{ on } \Gamma,\, v\in H^1(\Omega)\}.
 \end{align}
 Then, the weak formulation of (\ref{BVPBd1} )-(\ref{Recurrence2}) consists of 
 finding $(u,F_0,G_0\dots F_{L-1}, G_{L-1} )\in~\mathscr{S}$ such that the following equations are satisfied:
 \begin{enumerate}
\item Weak form of the governing equation
\begin{align}
a(u,v) - \sum_{l=0}^{L-1}c_l(F_l,v) - \sum_{l=0}^{L-1}d_l(G_l,v) =0, \label{WF1}\quad \mbox{for all $v\in \mathscr{S}_0$ }
 \end{align}
where
 \begin{align}
 & a(u,v)=\int_{\Omega} \left(\nabla u\cdot \nabla v - k^2 u v\right)\,d\Omega, \nonumber\\
 &c_l(F_l,v) = \mathscr{A}_l(kR)(F_l,v)_{S_R}, \quad  d_l(G_l,v) = \mathscr{B}_l(kR)(G_l,v)_{S_R}=0, \nonumber\\
&(F_l,v)_{S_R} = \int_{S_R} F_l v \,ds, \quad (G_l,v)_{S_R} = \int_{S_R} G_l v \,ds,\nonumber\\
&\mathscr{A}_l(kR)=-\frac{kH_1(kR)}{(kR)^l} - \frac{klH_0(kR)}{(kR)^{l+1}}\quad\mbox{and}\quad 
\mathscr{B}_l(kR)=-\frac{k(l+1)H_1(kR)}{(kR)^{l+1}} + \frac{kH_0(kR)}{(kR)^{l}},\nonumber
 \end{align}
 for $l=0,\dots L-1$. In the derivation of equation (\ref{WF1}), integration by parts after multiplication of the governing equation (\ref{BVPsc1}) by the test function $v$, and the continuity of the first radial derivative (\ref{BVPBd4}) at the artificial boundary $S_R$ have been used. 
 \item
Weak-form of the continuity of $u$  at $S_R$  (\ref{BVPBd3}) using Karp's expansion,
 \begin{align} 
 &  w_l(u,v_0)- \sum^{L-1}_{l=0}e_l(F_l,v_0) - \sum^{L-1}_{l=0}i_l(G_l,v_0)=0,\label{WF2} \quad \mbox{for all ${v}_0 $ in $H^0(S_R)$}
 \end{align}
  where
 \begin{align}
w_l(u,v_0)=(u,v_0)_{S_R}, \quad e_l(F_l,v_0) = \frac{H_0(kR)}{(kR)^l}(F_l,v_0)_{S_R},\quad  i_l(F_l,v_0) = \frac{H_1(kR)}{(kR)^l}(G_l,v_0)_{S_R}, \nonumber
\end{align}
for $l=0,\dots L-1$.

\item
A third equation may be obtained from the weak-form of the continuity of the second derivative at the artificial boundary (\ref{BVPBd5}). However, noticing that continuity of the Helmholtz operator at $S_R$ is also verified, then it is possible to replace condition (\ref{BVPBd5}) by the more convenient condition for the general finite element method (FEM) given by
\begin{eqnarray}
&&\left(\partial^2_{r} + \frac{1}{r}\partial _r+  \frac{1}{r^2}\partial^2_{\theta} + k^2\right)
 \left(H_0(kr)\sum_{l=0}^{L-1} \frac{F_l(\theta)}{(kr)^l} + H_1(kr)\sum_{l=0}^{L-1} \frac{G_l(\theta)}{(kr)^l}  \right)
 \bigg|_{r=R}  =\nonumber \\
 &&\qquad \left( \Delta_{r\theta} u + k^2 u \right) \bigg|_{r=R} = 0,\label{ContHelm}
 \end{eqnarray}
 The weak-form of (\ref{ContHelm}) is given by
 \begin{align}
  &\sum_{l=0}^{L-1} p_l(F_l,{\hat v})+ \sum_{l=0}^{L-1} q_l(G_l,{\hat v}) =0,\label{WF3} \quad \mbox{for all ${\hat{v}} $ in $H^1(S_R)$}
  \end{align}
where 
\begin{align}
&p_l(F_l,{\hat v}) =
\mathscr{P}_l(kR)(F_l,{\hat v})_{S_R} -  \frac{H_0(kR)}{R^2(kR)^l}b(F_l,{\hat v})\nonumber\\
 &q_l(G_l,{\hat v}) = \mathscr{Q}_l(kR)(G_l,{\hat v})_{S_R} -  \frac{H_1(kR)}{R^2(kR)^l}b(G_l,{\hat v}) \nonumber\\
&\mathscr{P}_l(kR)=\mathscr{E}_l(kR) + \frac{1}{R}\mathscr{A}_l(kR) +\frac{k^2H_0(kR)}{(kR)^l},\nonumber\\
&\mathscr{Q}_l(kR)=\mathscr{I}_l(kR) + \frac{1}{R}\mathscr{B}_l(kR) +  \frac{k^2H_1(kR)}{(kR)^l},\nonumber\\
&\mathscr{E}_l(kR)=-k^2\left[ \left(
\frac{1}{(kR)^l} - \frac{l(l+1)}{(kR)^{l+2}}\right)H_0(kR)  - \frac{2l+1}{(kR)^{l+1}} H_1(kR)\right],\nonumber\\
&\mathscr{I}_l(kR)=-k^2 \left[
\frac{2l+1}{(kR)^{l+1}}H_0(kR) +\left(\frac{1}{(kR)^l }- \frac{(l+1)(l+2)}{(kR)^{l+2}} \right)H_1(kR)
\right], \quad
   \mbox{for }  l=0,\dots L-1, \,\mbox{and} \nonumber\\
& b(w,v) = R^2\int_{S_R} w'\,v' \,ds.\nonumber
\end{align}
 The weak form (\ref{WF3}) is obtained by applying the radial derivatives of Helmholtz operator (\ref{ContHelm})  to Karp expansion, multiplying by a test function ${\hat v}\in H^1(S_R)$, integrating on $S_R$ and applying integration by parts with respect to the angular variable $\theta$.  
 \item
Weak form of the recurrence formulas (\ref{Recurrence2}) and (\ref{Recurrence1})
\begin{align}
&x_l(F_l, {\hat v}) + y_l(G_{l-1}, {\hat v}) = 0\label{WF4}, \quad \mbox{for all ${\hat{v}} $ in $H^1(S_R)$},\\
&r_l(G_l, {\hat v}) + t_l(F_{l-1}, {\hat v}) = 0\label{WF5}, \quad \mbox{for all ${\hat{v}}$ in $H^1(S_R)$}
\end{align}
where
\begin{align}
&x_l(F_l, {\hat v}) = 2l(F_l, {\hat v})_{S_R}, \quad y_l(G_{l-1}, {\hat v}) =l^2 ( G_{l-1}, {\hat v})_{S_R}  
  +b(G_{l-1}, {\hat v}),\nonumber\\
&r_l(G_l, {\hat v}) = 2l(G_l, {\hat v})_{S_R}, \quad t_l(F_{l-1}, {\hat v}) = - (l-1)^2 ( F_{l-1}, {\hat v})_{S_R},  
  +b(F_{l-1}, {\hat v})\nonumber
\end{align}
for $l=1\dots L-1.$
The above equations (\ref{WF1}),(\ref{WF2}), (\ref{WF3})-(\ref{WF5}) constitute the weak form of the BVP (\ref{BVPBd1})-(\ref{Recurrence2}). They can be used to simultaneously solve for the scattered field $u$, and the new families of unknowns, $F_l$ and, $G_l$ ($l=0,1,\dots L-1$), of the Karp's expansion defined at the artificial boundary.
\end{enumerate}

For the Neumann boundary condition at the obstacle bounding curve $\Gamma$ ($Z=1$, hard obstacle), $\mathscr{S}= H^1(\Omega)\times \overbrace{H^1(S_R)\times \dots \times H^1(S_R)}^{2L \rm\; times}$. The weak form of the governing equation changes to 
\begin{align}
a(u,v)
- \sum_{l=0}^{L-1} c_l(F_l,v)- \sum_{l=0}^{L-1} d_l(G_l,v)= - (\partial_{\bf n}\uinc ,v)_{\Gamma}, 
\quad \mbox{for all ${{v}} $ in $H^1(\Omega)$}\label{WFN1}
 \end{align}
 This weak formulation for the Neumann problem is completed with equations (\ref{WF2}), and (\ref{WF3})-(\ref{WF5})  that remain unchanged with respect to the Dirichlet case.

Finite element approximations of $u$, $F_l$, and $G_l$ are obtained by choosing finite-dimensional subspaces $\mathscr{S}^h$ of $\mathscr{S}$ and $\mathscr{S}_0^h$ of $\mathscr{S}_0$ with their respective bases $\{\phi_1,\dots \phi_{n_h}\}$ and $\{\psi_1,\dots \psi_{m_h}\}$.
By constructing finite elements to cover the physical domain $\Omega$, discretizing the weak forms (\ref{WF1}),(\ref{WF2}), (\ref{WF3})-(\ref{WF5}),  and using the above bases functions, we arrive to the following linear system (Dirichlet problem) in generic form:
$$A^h\vect{u}^h= \vect{b}^h$$
where
\begin{equation}
A^h= \begin{bmatrix} 
\vect{A}  & &\dots -\vect{C}_{l-1} & - \vect{D}_{l-1} &-\vect{C}_l & - \vect{D}_l\dots 
& -\vect{C}_{L-2} & - \vect{D}_{L-2}
& -\vect{C}_{L-1} & - \vect{D}_{L-1} \\
&\vect{W} \dots& \quad\,-\vect{E}_{l-1} & - \vect{I}_{l-1} &-\vect{E}_l & - \vect{I}_l\dots &&&-\vect{E}_{L-1} & - \vect{I}_{L-1} \\
 & &\vect{P}_{l-1} &  \vect{Q}_{l-1} &\vect{P}_l & \vect{Q}_l\dots &\vect{P}_{L-2} &  \vect{Q}_{L-2}
 &\vect{P}_{L-1} &  \vect{Q}_{L-1}  \\
& & &\vect{Y}_{l-1} &\vect{X}_l  && \vect{Y}_{L-2}&\vect{X}_{L-1}  \\
 &  &\vect{T}_{l-1} &  & & \vect{R}_l\dots &\vect{T}_{L-2} &&&  \vect{R}_{L-1} 
\end{bmatrix} 
\end{equation}

\setcounter{MaxMatrixCols}{20}
\begin{equation}
\vect{u}^h = \begin{bmatrix}
\vect{u}^h_{\Omega} &
\vect{u}^h_{S_R} &
\dots&
\vect{F}^h_{l-1}&
\vect{G}^h_{l-1}&
\vect{F}^h_{l}&
\vect{G}^h_{l}&
\dots &
\vect{F}^h_{L-2}&
\vect{G}^h_{L-2}&
\vect{F}^h_{L-1}&
\vect{G}^h_{L-1}
\end{bmatrix}^\top
\end{equation}
\begin{equation}
\qquad\qquad\vect{{b}}^h =\begin{bmatrix}
\vect{b}^h_{inc}&
0&
\dots&
0&
0&
0&
0&
\dots&
0&
0&
0&
0
\end{bmatrix}^\top
\end{equation}
As usual in finite elements, the entries in the block submatrices forming $A^h$ are obtained from the bilinear forms defining the weak forms (\ref{WF1}),(\ref{WF2}), (\ref{WF3})-(\ref{WF5}) acting on the basis functions. The uppercase letters ($\vect{A},\, \vect{C},\, \dots$) correspond to the lowercase letters found in the weak forms. The entries of the unknown vector ${\vect u}^h$ correspond to the values on the grid points of the unknown functions $u$, $F_l$ and $G_l$. Also, the only nonzeros entries $\vect{b}^h_{inc}$ of the vector $\vect{b}^h$ are obtained from the incident wave $u_{inc}$.


\section{The scattering problem in three dimensions. Wilcox's expansion}
\label{3DProblem}

The three dimensional scattering problem using Wilcox farfield expansion absorbing boundary condition (WFE) with L terms  was also introduced in \cite{JCP2017}. The corresponding equations in spherical coordinates are:
\begin{align}
& \Delta u + k^2 u = 0, \quad\qquad  &\text{in $\Omega$}, \label{BVP3D1} \\
& Z\, {\partial u}_n u + (1-Z) \,u = -\Big(Z
\,\partial_n u_{inc} + (1-Z) \,u_{inc}\Big) &\text{on $\Gamma$,}   \label{BVP3D2} \\
& u(R,\theta,\phi)=\frac{e^{ikR}}{kR}\sum_{l=0}^{L-1} \frac{F_l(\theta,\phi)}{(kR)^l}\label{BVP3D3} \\
&\partial_{r} u(R,\theta,\phi) = \frac{e^{i k R}}{k R} \sum_{l=0}^{L-1} \left( ik - \frac{l+1}{R} \right) \frac{F_{l}(\theta,\phi)}{(kR)^l}, \label{BVP3D4}\\
& 2ilF_l(\theta,\phi) = l(l-1) F_{l-1}(\theta,\phi) + \Delta_{\Sph} F_{l-1}(\theta,\phi), &\qquad l \geq 1, \label{BVP3D5}
\end{align}
where $\Delta_{\Sph}$ is the Laplace-Beltrami operator in the angular coordinates $\theta$ and $\phi$. See \cite{Bayliss01}.
Notice, that the WFE only has one unknown family of angular functions $F_l$ ($l=1\dots L-1$). As a consequence, only one  recurrence formula is needed. This recurrence formula and the continuity of the first radial derivative at the artificial boundary $S_R$ are sufficient to complete the equations defining the WFE-BVP. 

\subsection{Weak formulation and finite element approximation.}
\label{Weakformulation3D}
The weak form of (\ref{BVP3D1} )-(\ref{BVP3D5} ) for the Dirichlet problem ($Z=0$) can be obtained by first defining the function spaces 
\begin{align}
&\mathscr{S}=\{(u,F_0,\dots F_{L-1} ) \,|\,u=-u_{inc}  \mbox{ on } \Gamma,\, u\in H^1(\Omega),\, \nonumber
F_l\in H^1(S_R),\, l=0\dots L-2, \nonumber\\
&\qquad\,\,\,F_{L-1} \in H^0(S_R)\}\nonumber\\
&\mathscr{S}_0=\{v \,|\,v=0  \mbox{ on } \Gamma,\, v\in H^1(\Omega)\}.\nonumber
 \end{align}
 Then, the weak formulation consists of 
 finding $(u,F_0,\dots F_{L-1})\in~\mathscr{S}$ such that the following equations are satisfied:
 \begin{enumerate}
\item Weak form of the governing equation
\begin{align}
a(u,v) - \sum_{l=0}^{L-1}c_l(F_l,v) =0, \label{WF3D1}\quad \mbox{for all $v\in \mathscr{S}_0$ }
 \end{align}
where
 \begin{align}
 & a(u,v)=\int_{\Omega} \left(\nabla u\cdot \nabla v - k^2 u v\right)\,d\Omega, \quad
 c_l(F_l,v) = \frac{e^{ikR}}{(kR)^{l+1}}\left( ik - \frac{l+1}{R} \right) 
 (F_l,v)_{S_R},\nonumber\\
&(F_l,v)_{S_R} = \int_{S_R} F_l v \,ds,\qquad \mbox{for $l=0,\dots L-1$.}\nonumber
 \end{align}
 \item
Weak-form of the continuity of $u$  at $S_R$  (\ref{BVP3D3}) using Wilcox's expansion,
 \begin{align} 
 &  w_l(u,v_0)- \sum^{L-1}_{l=0}e_l(F_l,v_0),\label{WF3D2} \quad \mbox{for all ${v}_0 $ in $H^0(S_R)$}
 \end{align}
  where
 \begin{align}
w_l(u,v_0)=(u,v_0)_{S_R}, \quad e_l(F_l,v_0) = \frac{e^{ikR}}{(kR)^{l+1}}(F_l,v_0)_{S_R}, \qquad \mbox{for $l=0,\dots L-1$}.
\nonumber
\end{align}
\item
Weak form of the recurrence formula (\ref{BVP3D5}) 
\begin{align}
x_l(F_l, {\hat v}) + y_l(F_{l-1}, {\hat v}) = 0\label{WF3D3}, \quad \mbox{for all ${\hat{v}} $ in $H^1(S_R)$},
\end{align}
where
\begin{align}
&x_l(F_l, {\hat v}) = 2il(F_l, {\hat v})_{S_R}, \nonumber\\
&y_l(F_{l-1}, {\hat v}) = -l(l-1)( F_{l-1}, {\hat v})_{S_R}  + {R^2} \int_{S_R} \nabla_SF_{l-1}\cdot \nabla_S {\hat v}\,ds
\nonumber
\end{align}
for $l=1\dots L-1.$
The symbol $\nabla_{S}$ represent the gradient in the geometry of the sphere $S$.
The above equations (\ref{WF3D1})-(\ref{WF3D3}) constitute the weak form of the BVP (\ref{BVP3D1})-(\ref{BVP3D5}). They can be used to simultaneously solve for the scattered field $u$, and the new family of unknowns, $F_l$ ($l=0,1,\dots L-1$), of the Wilcox's expansion defined at the artificial boundary.
\end{enumerate}

For the Neumann boundary condition at the obstacle bounding curve $\Gamma$ ($Z=1$, hard obstacle), we define $\mathscr{S}= H^1(\Omega)\times \overbrace{H^1(S_R)\times \dots \times H^1(S_R)}^{L \rm\; times}$ and the weak form corresponding to the governing equation changes to 
\begin{align}
a(u,v)
- \sum_{l=0}^{L-1} c_l(F_l,v)= - (\partial_{\bf n}\uinc ,v)_{\Gamma}, 
\quad \mbox{for all ${{v}} $ in $H^1(\Omega)$}\label{WFN13D}
 \end{align}
 This weak formulation for the Neumann problem is completed with equations  (\ref{WF3D2})-(\ref{WF3D3})  that remain unchanged with respect to the Dirichlet case.

Finite element approximations of $u$ and $F_l$ are obtained by choosing finite-dimensional subspaces $\mathscr{S}^h$ of $\mathscr{S}$ and $\mathscr{S}_0^h$ of $H^1(\Omega)$ with their respective bases $\{\phi_1,\dots \phi_{n_h}\}$ and $\{\psi_1,\dots \psi_{m_h}\}$.
By constructing finite elements to cover the physical domain $\Omega$, discretizing the weak forms  (\ref{WF3D1})-(\ref{WF3D3}),  and using the above basis functions, we arrive to the following linear system (Dirichlet problem) in generic form:
$$A^h\vect{u}^h= \vect{b}^h$$
where

\begin{equation}
A^h= \begin{bmatrix} 
\vect{A}  & &\ -\vect{C}_{0} &- \vect{C}_{1} &\dots -\vect{C}_{l-1} & - \vect{C}_l\dots 
& -\vect{C}_{L-2}  
& -\vect{C}_{L-1} \\
&\vect{W} &-\vect{E}_{0} &-\vect{E}_{1}\dots &-\vect{E}_{l-1} & -\vect{E}_{l} \dots& - \vect{E}_{L-2}
& -\vect{E}_{L-1} \\
& & \vect{Y}_{0} &\vect{X}_1\dots  & \vect{Y}_{l-1}&\vect{X}_{l}  & \vect{Y}_{L-2}&\vect{X}_{L-1}  
\end{bmatrix} \label{stiff3D}
\end{equation}

\begin{equation}
\vect{u}^h =
\begin{bmatrix}
\vect{u}^h_{\Omega} &
\vect{u}^h_{S_R}& 
\vect{F}^h_{0}&
\vect{F}^h_{1} &
\dots &
\vect{F}^h_{l-1}&
\vect{F}^h_{l} &
\dots &
\vect{F}^h_{L-2} &
\vect{F}^h_{L-1}
\end{bmatrix}^\top
\end{equation}

\begin{equation}
\vect{{b}}^h =\begin{bmatrix}
\vect{b}^h_{inc}&
0&0&0&\dots &
0& 0&\dots & 0& 0
\end{bmatrix}^\top
\end{equation}
The description of the entries of the matrix $A^h$, the vector $\vect{u}^h$, and $\vect{ b}^h$ is completely analogous to their description in the 2D case.

We adopt IGA  to obtain the numerical solution of the finite element approximations obtained in this and the previous section in the artificially truncated computational domain. 
 This includes the numerical solutions for the two families of unknown angular functions ($F_l(\theta)$ and $G_l(\theta)$, $l=1\dots L-1$) defined on the artificial boundary by employing the same basis functions used to approximate both the geometry and the solution in the interior of the computational domain. Hence, the proposed methods is truly an isogeometric one providing a convenient platform to perform high order analysis and refinement. More precisely, both the domain geometry and the numerical solution in the interior and on the boundary are approximated using B-spline/NURBS basis functions.  More details about the NURBS basis function are found in the Appendix A. We denote this technique  by {\it IGA-FEABC}. In the following sections, we study both the convergence rate, and the accuracy of IGA-FEABC with {\it p}- and {\it h}- refinement and their dependence on the number of terms, $NT$, of the farfield expansions. 


\section{Numerical Experiments}
\label{sec:num}
We verify the high accuracy and the high order of convergence of the proposed IGA-FEABC by performing experiments in 2D and 3D acoustic scattering problems. They include plane wave scattering from an infinite circular cylinder, acoustic scattering from a prototype submarine in 2D, and scattering from a spherical scatterer. For the circular cylinder and spherical scatterer (axisymmetric case) comparison against the exact solutions allows to obtain the order of convergence and also a measure of the accuracy of the numerical method. We define the discretization density as the number of control points per wavelength and denote it with the symbol $n_{\lambda}$. We observe the dependence of the accuracy and the order of convergence on the number of terms $NT$ of the farfield expansions, the order $p$ of the basis functions employed, the discretization density $n_{\lambda}$, and the
number of degrees of freedom (DOF). We perform experiments for rather high frequencies such
as $k = 100$, and 350, very low frequencies such as $k = 0.01$, and for artificial boundaries located
as close as 0.05 distance from the scatterer boundary. For most of these experiments, we 
obtain highly accurate numerical solutions by appropriately adjusting $p$, $n_\lambda$, and
$NT$.


\subsection{Acoustic scattering from a circular cylindrical scatterer at middle and low frequencies}
\label{Circular}
First, we study the scattering of a plane wave propagating in the positive direction of the $x$-axis from a soft (Dirichlet BC) or hard (Neumann boundary BC) circular cylindrical scatterer  of radius $r_0=1$ for which the artificial boundary is at $R=2$. These benchmark problems have exact solutions in terms of eigenfunction expansions \cite{MartinBook}. As described in Section \ref{Weakformulation}, they are modeled by the equations (\ref{WF1}),(\ref{WF2}), (\ref{WF3})-(\ref{WFN1}) in weak-form. In this section we study the properties of the numerical solutions obtained from the application of the IGA-KFE technique to these equations.

\subsubsection {Accuracy and order of convergence}
In the Tables \ref{table:1}- \ref{table:3}, we report the order of convergence of the numerical solution for a moderate frequency $k=2\pi$ and for different orders $p$ of the basis functions with $h$-refinement. More precisely, the discretization density $n_{\lambda}$ is gradually increased so that a consistent order of convergence is achieved. The number of terms $NT$ of the KFE  is adjusted to obtain the best possible order of convergence for the various values of $p$. 
 For the annular circular region $\Omega$, we construct a mesh ($N\times m$ control points) with step-sizes $\Delta r=\frac{R-r_0}{N-1}$, and $\Delta \theta=\frac{2\pi}{m}$ in the radial and the angular directions, respectively. Hence, the number of control points in the radial and angular directions are $N=(R-r_0)\nl$ and $m= 2\pi\nl$, respectively. For the convergence analysis, we also define $h=\Delta r=r_0\Delta \theta$.
\begin{table}[h!]
\centering
\begin{tabular}{||c c c c c||}
\hline
   $n_{\lambda}$    & Grid size & $h= r_{0} \Delta \theta = \Delta r$   & $L^2$-norm Rel. Error   & Observed order \\ [0.5ex]
   \hline\hline
    $20$   & $20\times 125$  & $0.05027$  & $2.95\times 10^{-5}$  & $ $    \\
 $23$   & $23\times 145$  & $0.04333$  & $1.86\times 10^{-5}$  & $ 3.12$    \\
 $25$   & $25\times 157$  & $0.04002$  & $1.47\times 10^{-5}$  & $ 2.93$    \\
 $29$   & $29\times 181$  & $0.03471$  & $9.65\times 10^{-6}$  & $ 2.96$    \\
$32$   & $32\times 201$  & $0.03126$  & $6.94\times 10^{-6}$  & $ 3.15$    \\
$    $  && & Least squares fit \,\, =    &3.03\\
[1ex]
\hline
\end{tabular}
\caption{Order of convergence of the numerical solution at the artificial boundary for order $p=2$  basis functions, using KFE with $NT=11$ terms}
\label{table:1}
\end{table}
\begin{table}[h!]
\centering
\begin{tabular}{||c c c c c||}
\hline
   $n_{\lambda}$    & Grid size & $h= r_{0} \Delta \theta = \Delta r$   & $L^2$-norm Rel. Error   & Observed order \\ [0.5ex]
   \hline\hline
    $28$   & $28\times 173$  & $0.03632$  & $1.35\times 10^{-6}$  & $ $    \\
 $30$   & $30\times 189$  & $0.03324$  & $9.33\times 10^{-7}$  & $ 4.17$    \\
 $32$   & $32\times 201$  & $0.03126$  & $7.25\times 10^{-7}$  & $ 4.10$    \\
 $34$   & $34\times 213$  & $0.02950$  & $5.73\times 10^{-7}$  & $ 4.06$    \\
$36$   & $36\times 225$  & $0.02793$  & $4.59\times 10^{-7}$  & $ 4.03$    \\
$    $  && & Least squares fit \,\, =    &4.04\\
[1ex]
\hline
\end{tabular}
\caption{Order of convergence of the numerical solution at the artificial boundary for order $p=3$  basis functions, using KFE with $NT=10$ terms}
\label{table:2}
\end{table}
\begin{table}[h!]
\centering
\begin{tabular}{||c c c c c||}
\hline
   $n_{\lambda}$    & Grid size & $h= r_{0} \Delta \theta = \Delta r$   & $L^2$-norm Rel. Error   & Observed order \\ [0.5ex]
   \hline\hline
    $26$   & $26\times 161$  & $0.03903$  & $6.14\times 10^{-8}$  & $ $    \\
 $28$   & $28\times 173$  & $0.03632$  & $4.15\times 10^{-8}$  & $ 5.45$    \\
 $30$   & $30\times 189$  & $0.03324$  & $2.59\times 10^{-8}$  & $5.34$    \\
 $32$   & $32\times 201$  & $0.03126$  & $1.88\times 10^{-8}$  & $ 5.16$    \\
$34$   & $34\times 213$  & $0.02950$  & $1.42\times 10^{-8}$  & $ 4.87$    \\
$    $  && & Least squares fit \,\, =    &5.25\\
[1ex]
\hline
\end{tabular}
\caption{Order of convergence of the numerical solution at the artificial boundary for order $p=4$ basis functions, using KFE with $NT=11$ terms}
\label{table:3}
\end{table}
\begin{figure}[h!]
\begin{center}
\includegraphics[width=0.5\linewidth]{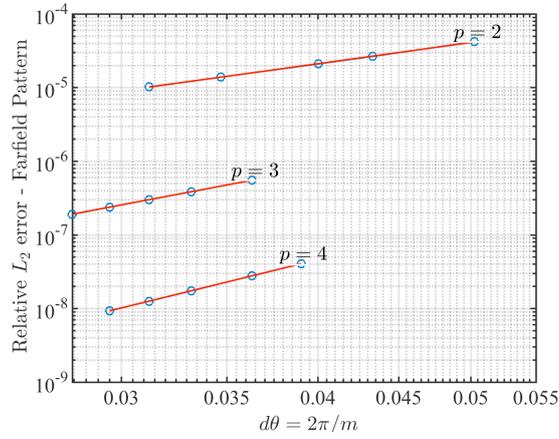} 
\end{center}
\caption{Least squares fitting lines for the observed order of convergence in Tables \ref{table:1} \ref{table:2} \ref{table:3}. }
\label{OrdConv}
\end{figure}

The Tables \ref{table:1}-\ref{table:3} reveal that employing bases of order $p$ in the IGA-KFE method leads to a numerical technique of $O(h^{p+1})$ under $h$-refinement, as expected. In each experiment, we needed to increase the number of terms $NT$ of the KFE until the order of convergence $p+1$ was achieved. Hence, it is possible to achieve the accuracy of a high order numerical method over the entire
computational domain including the artificial boundary by employing an appropriate number of
terms in the KFE used at the artificial boundary. A comparison of the least squares fit lines for the various experiments reported in Tables  \ref{table:1}-\ref{table:3} {are depicted} in Fig. \ref{OrdConv}. This figure shows the three least square lines with their respective slopes {illustrating} the corresponding order of convergence.

\begin{figure}[!h]
\begin{center}
\includegraphics[width=0.49\textwidth]{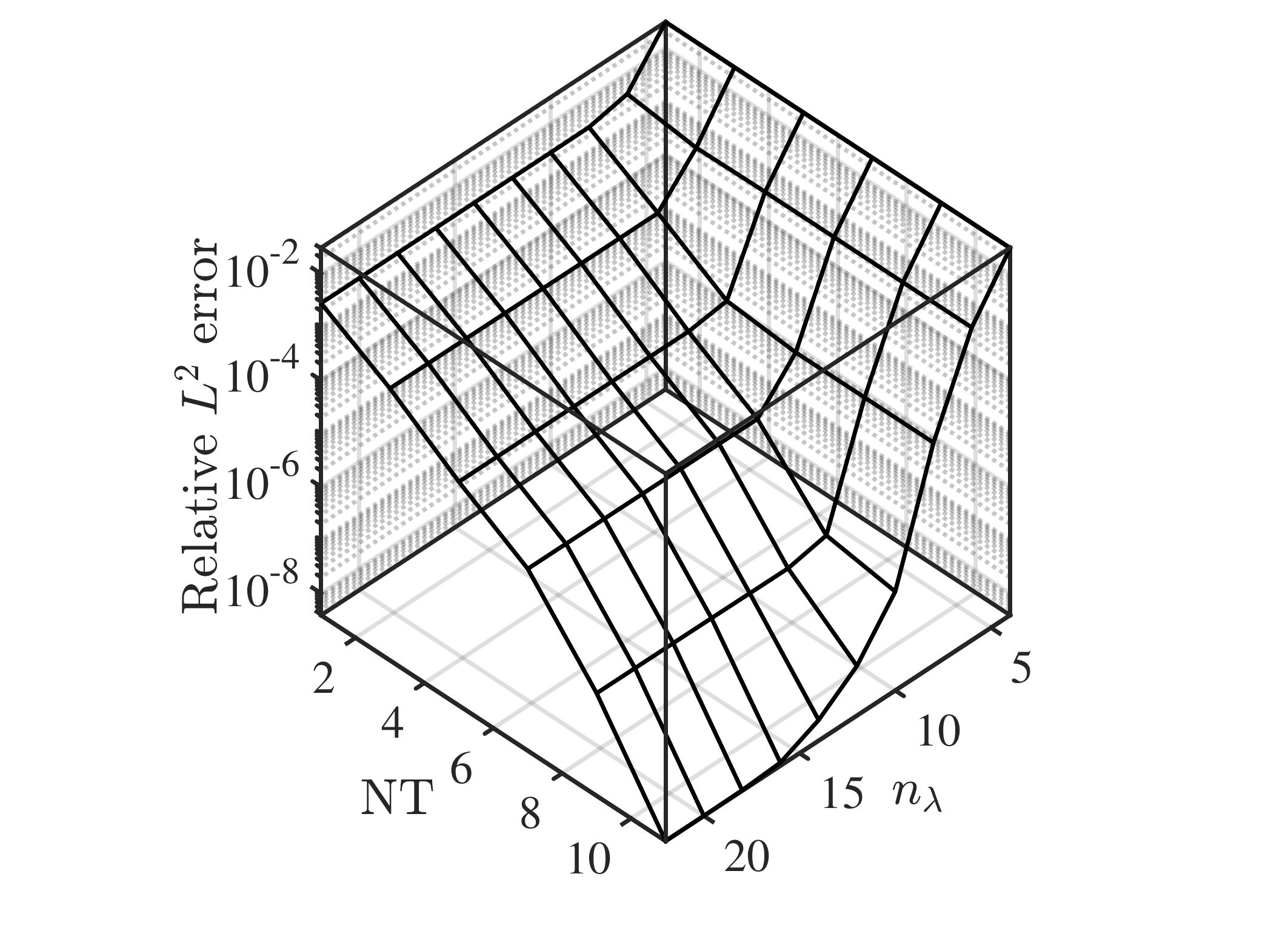} 
\includegraphics[width=0.49\textwidth]{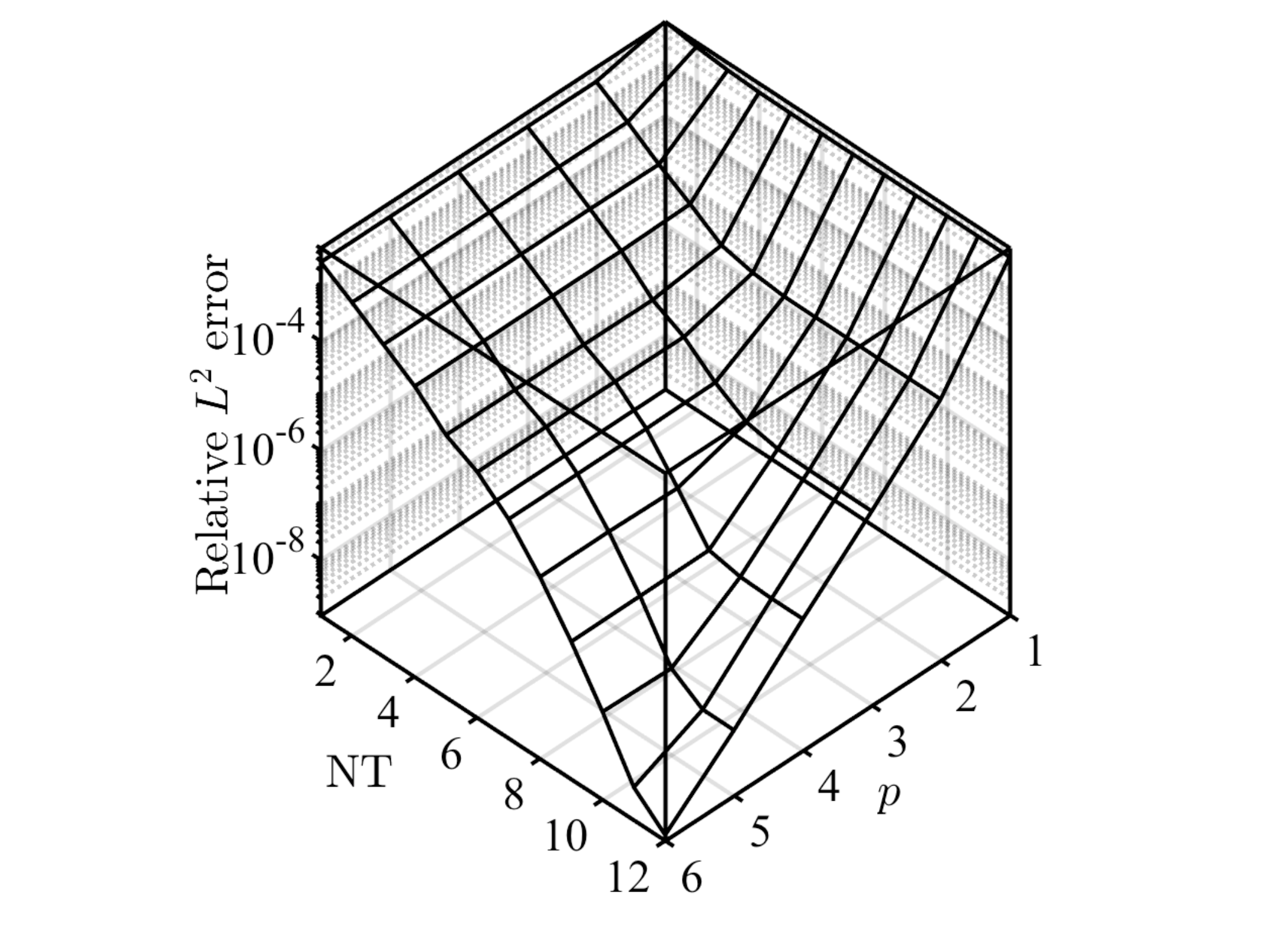} 
\caption{Farfield Pattern relative $L^2$ error for $p=6$ (left) and
relative $L^2$ error for $\nl=18$ (right).} 
\label{Surfp2D}
\end{center}
\end{figure}

In most of our experiments, we also compute numerical approximations of the {\it Farfield Pattern (FFP) of the scattered wave}. This is defined by the angular function present in the dominant term of the asymptotic expansion of the scattered wave when $r\rightarrow\infty$. For its analytical expression in 2D and 3D see \cite{MartinBook}. Also, its efficient calculation in 2D from the numerical solution of the scattered wave is well explained in \cite{JCP2017}. The dependence of the IGA-KFE accuracy on $n_{\lambda}$ and $NT$ is illustrated by the surface graph shown in Fig. \ref{Surfp2D} (left) for basis of order $p=6$, with $k=2\pi$ and $R=2$. We observe that the relative $L^2$ error in the computation of the FFP  decreases as $\nl$ increases. However, this occurs only up to certain $\nl$ value from which appreciable changes are not observed by increasing $\nl$. However, as $NT$ is further increased the error continue decreasing as $\nl$ increases. The minimum error shown in this figure is approximately $5\times 10^{-9}$ which corresponds to $\nl = 22$ and $NT =11$. 
  
 In the right side of Fig. \ref{Surfp2D}, a decaying of the L2 relative error is observed for $n_{\lambda}= 18$ fixed
when $p$ and $NT$ increase. However, there is not much error reduction for $NT\le 6$ even if $p$ is increased.
But, as $NT$ is made greater than $6$, the error decreases
  at a faster rate with $p$ refinement until it reaches a minimum value of approximately $10^{-8}$ for $p=6$ and $NT=12$. These results verify the unusual high accuracy that can be obtained by employing the IGA-KFE technique proposed in this work. In principle, further increase of $NT$ and $p$ would result in even smaller error. But, it might be necessary to employ an iterative solver to solve the resulting linear system.
  
 \subsubsection{Comparison of IGA-KFE with other numerical techniques and absorbing boundary conditions}
  In this section, we report on a series of experiments to highlight the advantages of the KFE and the IGA-KFE over some well-known ABCs and similar numerical methods, respectively. In the performed experiments, we consider a sound soft scatterer and compute the numerical solution in the region enclosed by the scatterer and the artificial boundary located at $R=2$, for a frequency $k=10$.
    
\begin{figure}[!h]
\begin{center}
\includegraphics[width=0.49\textwidth]{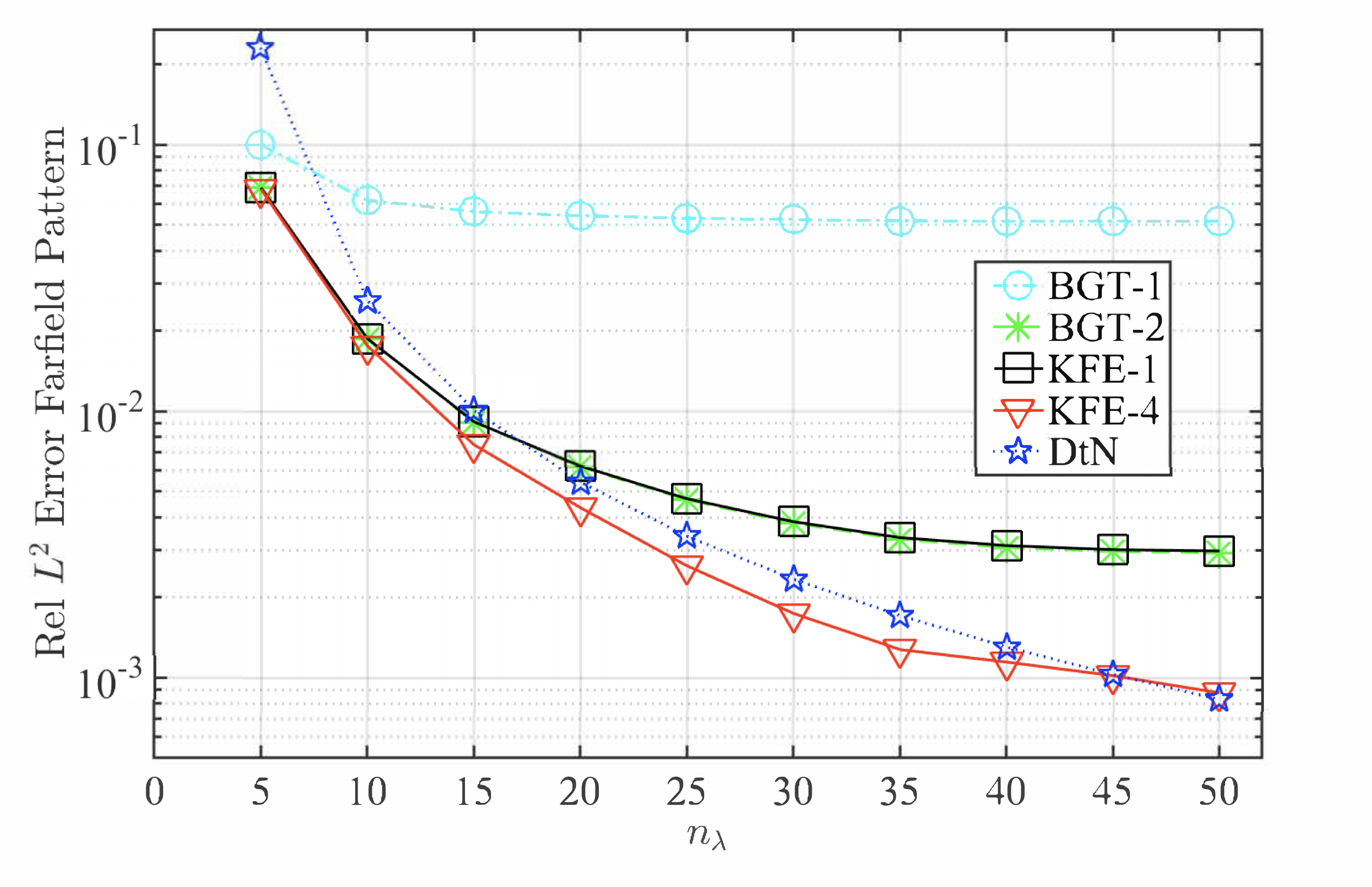}  
\includegraphics[width=0.49\textwidth]{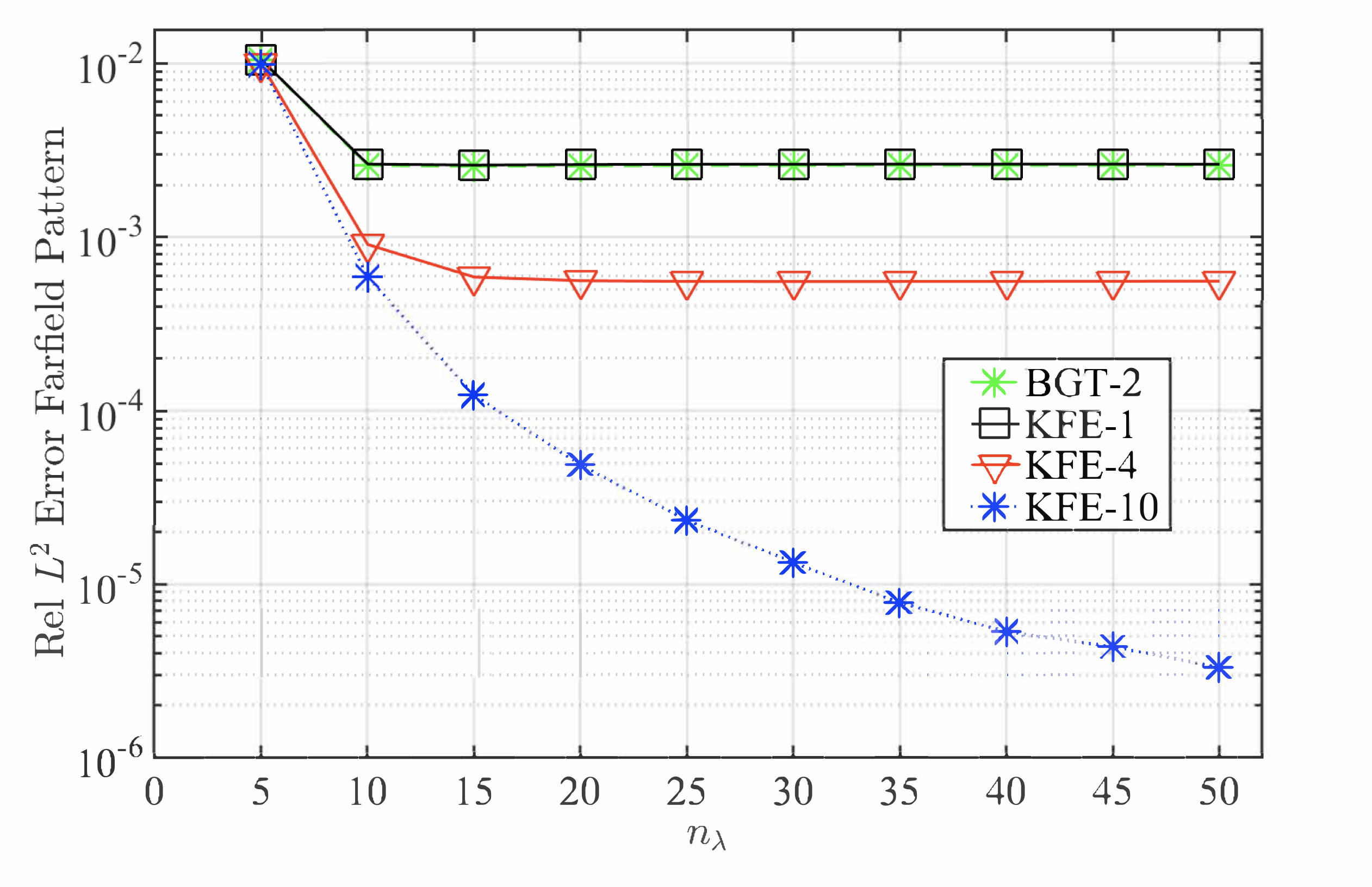} 
\caption{Convergence comparison for   $k=10$, and $R=2$ of IGA combined with various ABC and FD$_2$-DtN. The graph at the left is  for NURBS bases of order $p=1$  and the one to the right is for  NURBS bases of order $p=2$ } 
\label{ABCComparison1}
\end{center}
\end{figure}

  First in Fig.\ref{ABCComparison1} (left side), we present the results of experiments combining IGA with the following absorbing boundary conditions: BGT-1, BGT-2, KFE-1, and KFE-4. The number next to BGT specifies the order of approximation to the Sommerfeld radiation condition while the one next to KFE represents the number of terms  $NT$  in Karp's expansions. The subindex $k$, in the notation IGA$_k$ for $k=1,2,\dots n$, is introduced to designate the order $p$ of the NURBS basis employed by the IGA technique.
We also include in Fig. \ref{ABCComparison1} the performance of a second order finite difference approximation 
  combined with a Dirichlet to Neumann absorbing boundary condition (FD$_2$-DtN). It is observed that IGA$_1$-BGT-1 reaches a stagnation 
  value for $n_{\lambda}\ge 15$. In an attempt
to correct this, we also tried with $p=2$. For low $n_{\lambda}$, the error is smaller, but it is also not decreasing
by refining. These experiments show that BGT-1 provides a very poor approximation at the
boundary. Therefore, the IGA$_{1,2}$-BGT-1 computation with the coarsest grid cannot be improved by
refining it ($h$-refinement) or increasing the order of the basis ($p$-refinement) because the boundary
error is dominating the overall error.

On the other hand, the approximations obtained for IGA$_1$-BGT-2 and IGA$_1$-KFE-1 are almost identical. This is expected because BGT-2 can be considered as an asymptotic version of KFE-1 when $R\rightarrow\infty$. It is noticeable that the approximation greatly improves compared with the previous computation, as shown in the graph, but the error changes little beyond $\nl \ge 30$. Again, the boundary error starts dominating the overall computation about this discretization density.
Finally, we employ $NT = 4$ terms in Karp's expansion, we observe smaller errors than IGA$_1$-BGT-2 and IGA$_1$-KFE-1. Additionally, there is not stagnation point for the range of $\nl$ in the figure. Actually for a sequence of $\nl =10,13,16,19, \mbox{and}\,\, 22$, IGA$_1$-KFE-4 exhibits quadratic convergence. This is the optimum results that can be obtained for an IGA$_1$ computation. Therefore, an increase in the number of terms beyond $NT=4$ for IGA$_1$ will not produce better results for this $\nl$ range. This is confirmed by applying a centered second order finite difference combined with an exact DtN absorbing boundary condition to this scattering problem. In fact, the error curves for both techniques have very similar behavior, as shown in Fig. \ref{ABCComparison1}.
  
To obtain lower errors for the same range of discretization density $[0,50]$, it is necessary to employ a technique with an order of convergence higher than 2.  A natural choice is to employ and IGA$_2$ method for the interior, i.e., the IGA method with a basis of order $p=2$. Because, this should reach a third order of convergence for the computation in the interior. 
 For this purpose, we combine the  IGA$_2$ method with the ABCs: BGT-2, KFE-1, KFE-4, and KFE-10. 
 In the right side of Fig. \ref{ABCComparison1}, we show the relative $L^2$ errors for the farfield pattern obtained from these combined methods. It is observed that the error produced at the artificial boundary by the ABCs: BGT-2, KFE-1, and KFE-4 dominates the computation beyond $\nl=10$ for BGT-2 and KFE-1 and for $\nl=15$ for KFE-4.  As a consequence, the error does not decrease for greater $\nl$. On the contrary for KFE-10, the relative $L^2$ error decreases until it reaches a minimum value of about $2\times10^{-5}$ when $\nl=50$. This means that the error due to the absorbing boundary condition KFE-10 at the artificial boundary is smaller than the one produces by the IGA$_2$ in the interior of the computational domain. Actually, the combined method IGA$_2$-KFE-10 has an order of convergence equal to 3 within this range of $\nl$.  
 
\begin{figure}[!h]
\begin{center}
\includegraphics[width=0.49\textwidth]{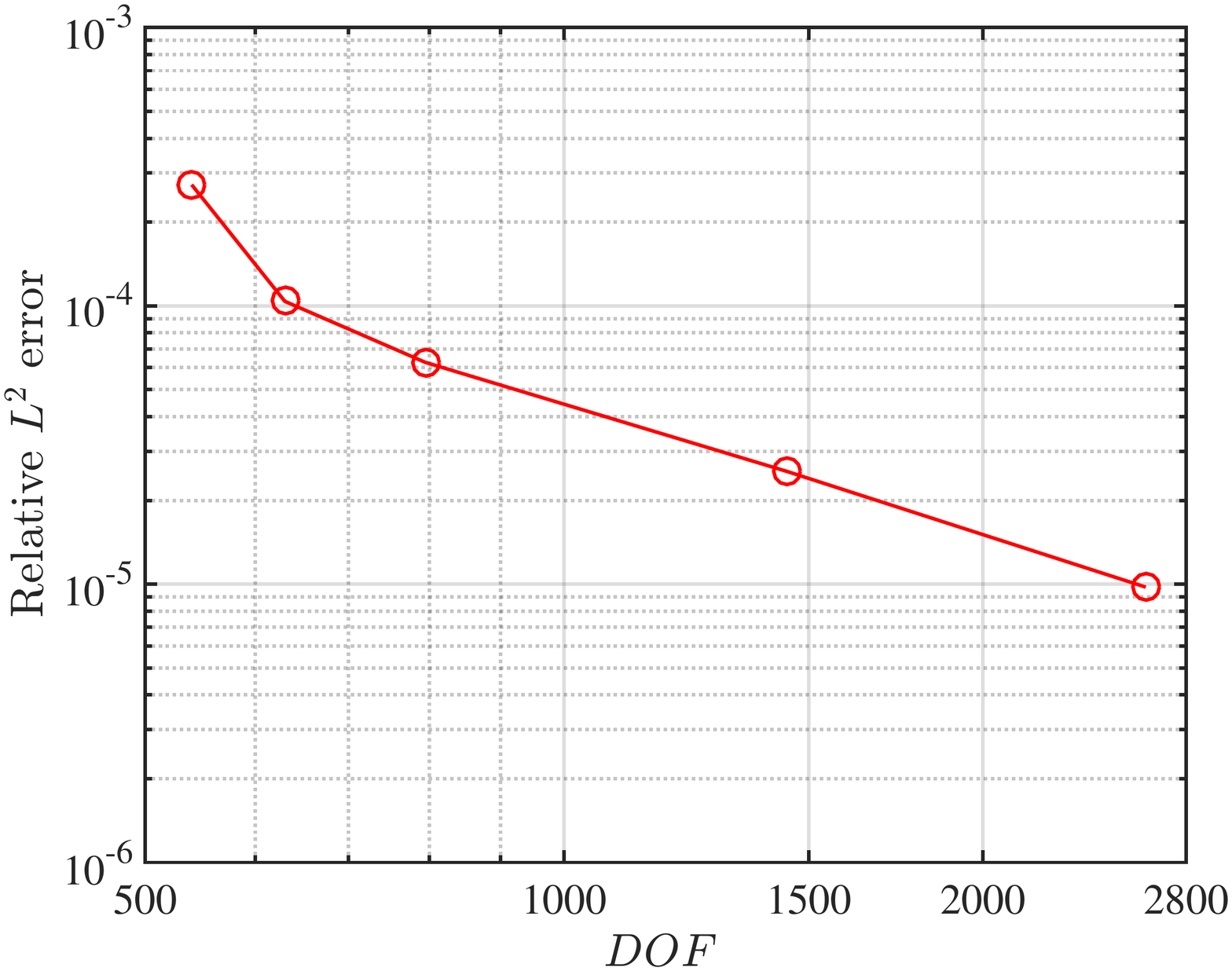}  
\includegraphics[width=0.49\textwidth]{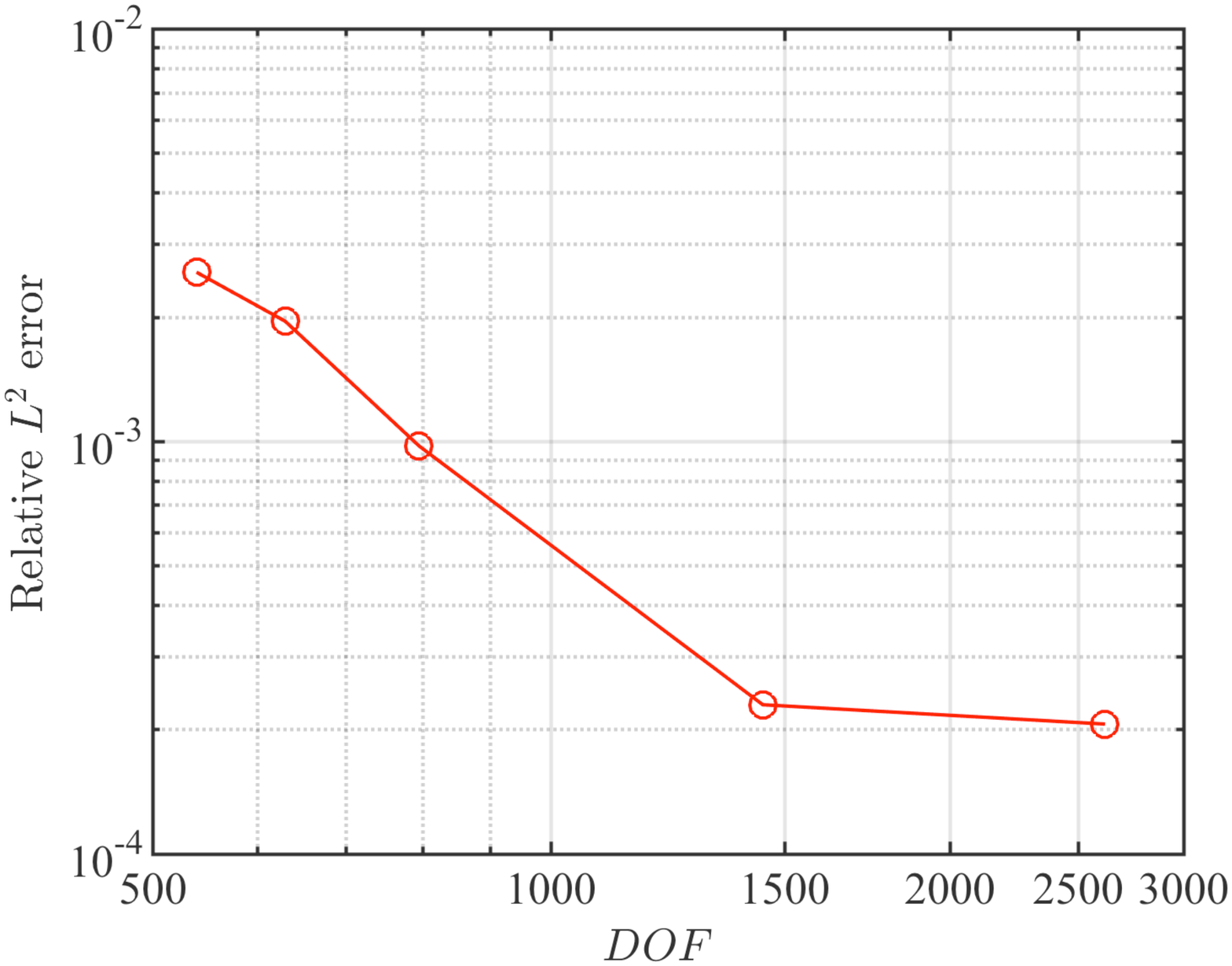} 
\caption{Convergence of IGA-KFE against DOF. Results for $k=1$ and $p=2$ (left), and for $k=10$ and $p=6$ (right).}
\label{DOFGraph1}
\end{center}
\end{figure}

 We also conducted another set of experiments employing IGA-KFE for the scattering from a sound-hard circular cylinder when $R=2$. We calculated the relative $L^2$ error in the computational domain against the number of degrees of freedom (DOF). In Fig. \ref{DOFGraph1} (left), the results are shown for a frequency $k=1$ with a basis of order $p=2$ and $NT=3$ terms of the KFE. Also in Fig. \ref{DOFGraph1} (right), we present the results for a frequency $k=10$ with $p=6$ and $NT=2-9$. 
 
In \cite{Antoine2009}, the authors sought the numerical solution for the same problem with the BGT-2 absorbing boundary condition located at $R=2$ using the Plane Wave enriched finite element method based on triangular quadratic finite elements (PWT6). The number of plane waves used was $q=1-4$. To avoid the error due to any ABC, they compared their numerical solution against the exact solution for the BVP with the BGT-2 absorbing boundary condition. The error obtained by applying this technique for $k=1$ was $4\times 10^{-4}$ for all the $DOF=840-3360$. In our experiments, we employed the KFE absorbing boundary condition and compared against the exact solution for the scattering problem. Although the error from the ABC were not avoided, we still obtained an error of $10^{-4}$ for $k=1$  for $DOF=630$. This error decreased to $9.8\times 10^{6}$ for $DOF=2619$,  as shown in Fig. \ref{DOFGraph1} (left).
For $k=10$, the results reported in \cite{Antoine2009} were $5\times 10^{-4}$ for $DOF=3360$ and $4\times 10^{-4}$ for $DOF=12480$. Our results using IGA$_6$ with $NT=2-9$ comparing against the exact solution were as low as $10^{-4}$ for $DOF=2387$, as illustrated in Fig. \ref{DOFGraph1} (right).
This error can be further decreased to $6\times 10^{-5}$ for $DOF=2759$, and even more to $4.12\times 10^{-6}$ for $DOF=3492$.

In \cite{Dinachandra-Raju2018}, a technique consisting of a Partition of Unity Isogeometric Analysis (PUIGA) coupled with BGT-2 at the artificial boundary $R=2$ was applied to the same 2D scattering problem with a frequency $k=10$. They employed $q=8$ plane waves to enrich their basis functions. The best relative $L^2$ error reported for the entire domain was $3.27\times 10^{-3}$ for $DOF=1280$ while the error obtained after applying IGA$_6$-KFE-2 was $2.6\times 10^{-3}$ for $DOF=1104$, as seen in Fig. \ref{DOFGraph1} (right). 
The error from the PUIGA coupled with BGT-2 will eventually reach a stagnation point such as the one shown in Fig. \ref{ABCComparison1} (right)  for the IGA technique coupled with BGT-2. This is due to the error generated at the artificial boundary by BGT-2. As a consequence even if the DOF is increased beyond certain value, the relative error will not decrease. 
These experiments show the advantage of the proposed method IGA-KFE over similar techniques which use low order ABC such as the BGT-2 at the artificial boundary. 

\subsubsection{Acoustic scattering from a circular cylinder at very low frequencies}
\label{CircularLow}
In this section, we report highly accurate results obtained by applying IGA-KFE method to the acoustic scattering at very low frequencies such as $k=0.01$.  We were inspired by similar work done by Grote and Keller \cite{Grote-Keller01} and Turkel et al. \cite{Turkel-Farhat2004}. These authors found
that employing BGT-2 as an ABC to obtain numerical solutions for the 2D acoustic scattering at very low frequencies leads to approximations several orders of magnitude different than the exact solution. This is due to the asymptotic character of the BGT-2 in two dimensions. They also found more accurate numerical solutions by employing as an ABC a second order differential operator, BGTH, which annihilates the leading order term of the Karp's expansion (\ref{BVPBd3}). They used for the discretization of the computational domain a second order finite difference and linear finite elements, respectively. Therefore, their results were limited to second order convergence at the most. 

In Table \ref{LowFreq}, 
we present the relative $L^2$ error over the entire domain obtained from the application of the IGA-BGT-2 and IGA-KFE technique to the 2D scattering from a sound-soft obstacle. The first column describe the radius $R$ identifying the location of the artificial boundary, the second column contains the number of elements in the radial direction, $N_e$ times the number of elements in the angular direction $m_e$.  The third column contains the relative $L^2$ errors when the IGA$_1$-BGT-~2 is applied and the remainder columns contains the relative $L^2$ error for the various combined methods employed depending of the $NT$ number of terms in Karp's expansion and the order $p$ of the bases. For comparison purpose, we performed experiments with the same data used in Table VII in \cite{Turkel-Farhat2004}. For IGA$_1$-BGT-2 and IGA$_1$-KFE-1, we obtained very similar results as those reported in Table VII for a linear FEM coupled with BGT-2 and BGTH, respectively. This is expected because IGA$_1$ and linear FEM have second order of convergence and BGTH is the differential operator which annihilates the two terms of Karp's expansion for $NT=1$. 
 \begin{table}[h!]
\centering
  \begin{adjustbox}{max width=\textwidth}
\begin{tabular}{c c c c c c c c}
\hline
   $R$    & $N_e\times m_e$ & IGA$_1$-BGT-2   & IGA$_1$-KFE-1 & IGA$_2$-KFE-1 & IGA$_2$-KFE-3 & IGA$_5$-KFE-3  & IGA$_{10}$-KFE-3 \\ [0.5ex]
   \hline
 $1.1$   & $10\times 60$  & $1.05\times 10^{-2}$  & $3.69\times 10^{-6}$  & $9.15\times 10^{-7}$  & $7.41\times 10^{-8}$  & $6.20\times 10^{-12}$ & $1.26\times 10^{-12}$\\
 $2$   & $10\times 60$  & $6.11\times 10^{-2}$  & $5.66\times 10^{-5}$  & $1.57\times 10^{-6}$  & $1.37\times 10^{-6}$  & $2.30\times 10^{-10}$ & $8.68\times 10^{-13}$\\
 $3$   & $20\times 60$  & $8.17\times 10^{-2}$  & $6.62\times 10^{-5}$  & $8.43\times 10^{-7}$  & $7.67\times 10^{-7}$  & $1.24\times 10^{-10}$ & $7.36\times 10^{-13}$\\
 $5$   & $40\times 60$  & $9.58\times 10^{-2}$  & $7.38\times 10^{-5}$  & $4.17\times 10^{-7}$  & $3.97\times 10^{-7}$  & $6.1\times 10^{-11}$ & $3.75\times 10^{-13}$\\
[1ex]
\hline
\end{tabular}
\end{adjustbox}
\caption{Relative $L^2$ error over the computational  domain for sound-soft acoustic scattering for a  frequency $k=0.01$. }
\label{LowFreq}
\end{table}

\begin{figure}[!h]
\begin{center}
\includegraphics[width=0.6\linewidth]{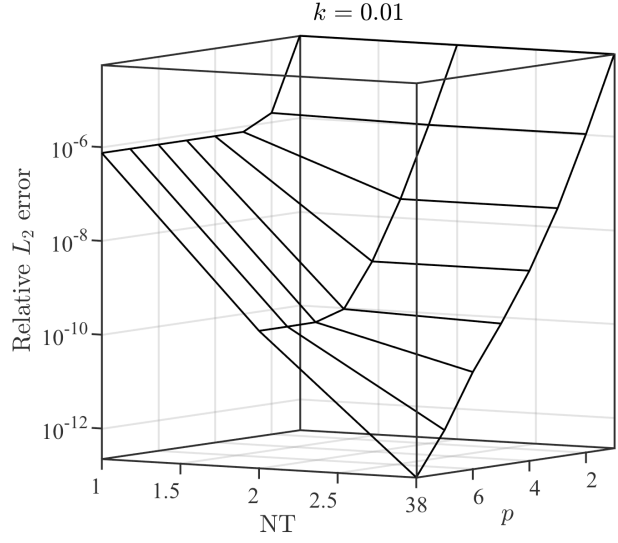}  
\caption{The entire domain relative $L^2$ error for the the sound-soft circular cylinder scatterer at a frequency $k=0.01$. }
\label{LowFreqGraph}
\end{center}
\end{figure}

 As it is shown in Table  \ref{LowFreq},  the error decreases almost to machine precision by implementing $p$-refinement. More precisely, the order of the IGA bases is increased from $p=1$ to $p=10$ while only $NT=3$ terms at  most are employed for the Karp's expansion at the artificial boundary. Therefore, the computational cost due to the use of the KFE absorbing boundary condition is minimal. In Fig. \ref{LowFreqGraph}, the evolution of the relative $L^2$ error is depicted in terms of $p$ and $NT$. The artificial boundary is located at $R=2$ and the number of elements is $10\times 60$.
 The remarkable reduction of the error to almost machine precision by $p$-refinement and by increasing $NT$ in the IGA-KFE technique is clearly evident from this graph.

\subsubsection{Acoustic scattering from a circular cylinder with a very close artificial boundary}
\label{CircularClose}
\begin{figure}[!h]
\begin{center}
\includegraphics[width=0.6\textwidth]{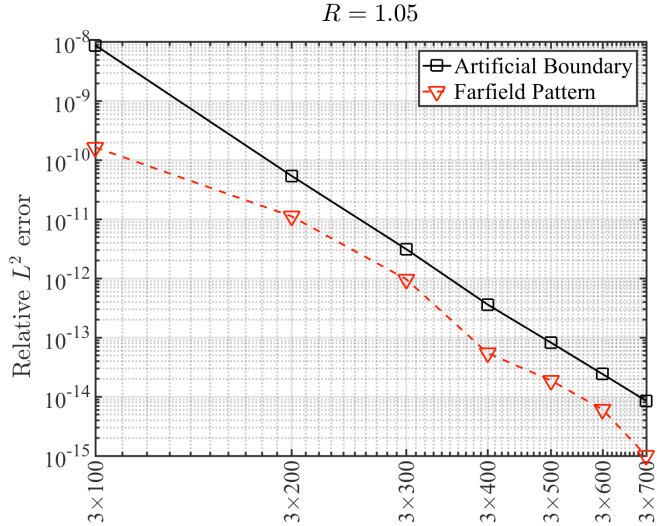}  
\caption{Relative $L^2$ error for a very close artificial boundary. }
\label{R105}
\end{center}
\end{figure}

 Another remarkable result showing the high accuracy of the IGA-FEABC technique is described in this section. In fact, we applied the combined method to the extreme problem where the artificial boundary radius ($R=1.05$) is chosen almost on top of the radius ($r_0=1$) of the circular scatterer. As a consequence, the domain of computation is very small which is an ideal situation to improve the  efficiency of the computational method. In the Fig. \ref{R105} two curves are graphed. One of them (solid line) corresponds to the 
 the evolution of the relative $L^2$ error at the artificial boundary under $h$-refinement. The other corresponds to the relative $L^2$ error for the farfield pattern (discontinuous line). 
 These two curves  were generated using only three elements 
 $N_e=3$, in the radial direction while in the angular direction the elements varied from $m_e=100$ to $m_e=700$. We also maintained fixed the number of karp's expansion terms as $NT=24$, and the degree of the NURBS basis employed was $p=6$  for a frequency $k=2\pi$ in these experiments.
 
 As it can be seen from the Fig. \ref{R105}, these relative errors are extraordinary low. For the farfield pattern, the error is already about $10^{-10}$ even for the coarsest mesh ($3\times100$).  Remarkably, for the finest mesh ($3\times 700$) the error, $9.99669\times 10^{-16}$, is at the level of machine precision for a computer using double precision. An appreciation of the magnitude of this result can be obtained by comparing it with the error, $1.304671\times 10^{-2}$,  obtained by applying the IGA technique coupled with the BGT-2 absorbing boundary condition to the same problem with identical input data.  It is worth noting that the BGT-~2 absorbing boundary condition was reported by Laghrouche \cite{Laghrouche2000} to give the higher accuracy among all the ABCs employed for similar  scattering problems. 
 
 In Table \ref{Time}, we report the CPU times spent in solving the linear systems corresponding to two experiments of the same scattering problem with artificial boundaries of radius $R=1.05$ and $R=5$, respectively.  The goal was to compare the time invested in each case to reach a relative $L^2$ error of order $10^{-8}$ for the computation of the scattered field at the absorbing boundary.
 The results showed that employing an artificial boundary of radius $R=1.05$ only required $7\%$  of the time employed by the artificial boundary of radius $R=5$ to reach the same precision. The other parameters $k=2\pi$, $p=5$ and $\nl = 20$ are the same for both experiments. These results were generated on a Surface Pro 4 computer with Intel� Core i5-6300 CPU @2.40 GHz 2.5 GHz
with 8 GB RAM. 
 This shows the monumental computational time savings obtained by the proposed combined method IGA-FEABC.
 
  \begin{table}[h!]
\centering
  \begin{adjustbox}{max width=\textwidth}
\begin{tabular}{c c c c c }
\hline
   $R$   & $DOF$ & $\mbox{Rel } L^2 \mbox{ AB}$  &  Time( sec) & NT \\ [0.5ex]
   \hline
 $1.05$   & $4500$  & $5.62215\times 10^{-8}$  & $0.244819$  & 15 \\
 $5$   & $11250$  & $8.16815\times 10^{-8}$  & $3.446884$  & 5 \\
[1ex]
\hline
\end{tabular}
\end{adjustbox}
\caption{Time comparison to reach a precision of $10^{-8}$ at the absorbing boundary for the IGA-KFE for artificial boundaries of radius $R=1.05$ and $R=5$. }
\label{Time}
\end{table}

\subsection{Plane wave scattering from a prototype submarine }
\label{Submarine}
We consider the scattering of a plane wave propagating in the positive direction of the $x$-axis from a sound-soft prototype two-dimensional submarine. 

\begin{figure}[!h]
\begin{center}
\includegraphics[width=0.47\linewidth]{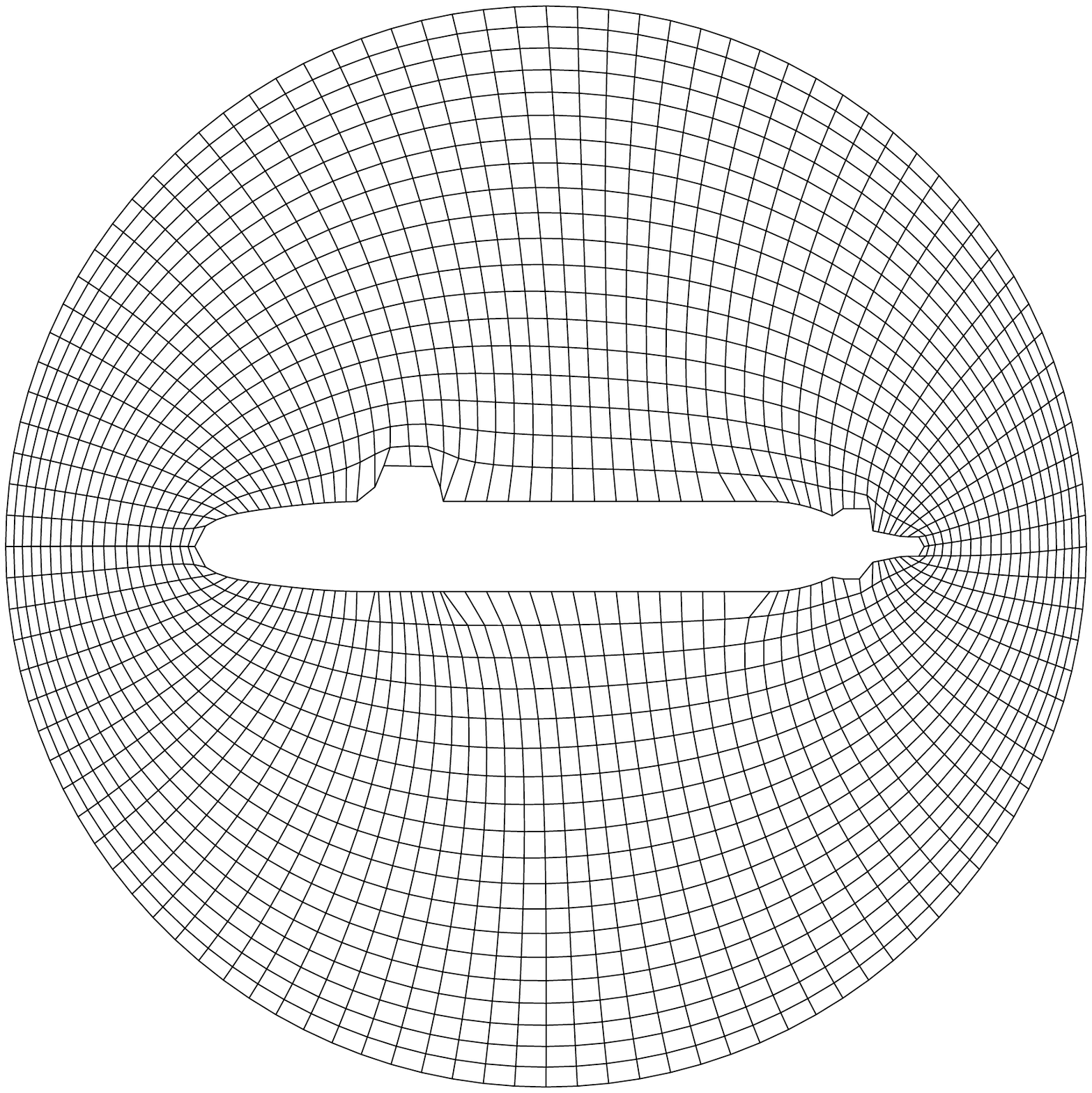}  
\includegraphics[width=0.49\linewidth]{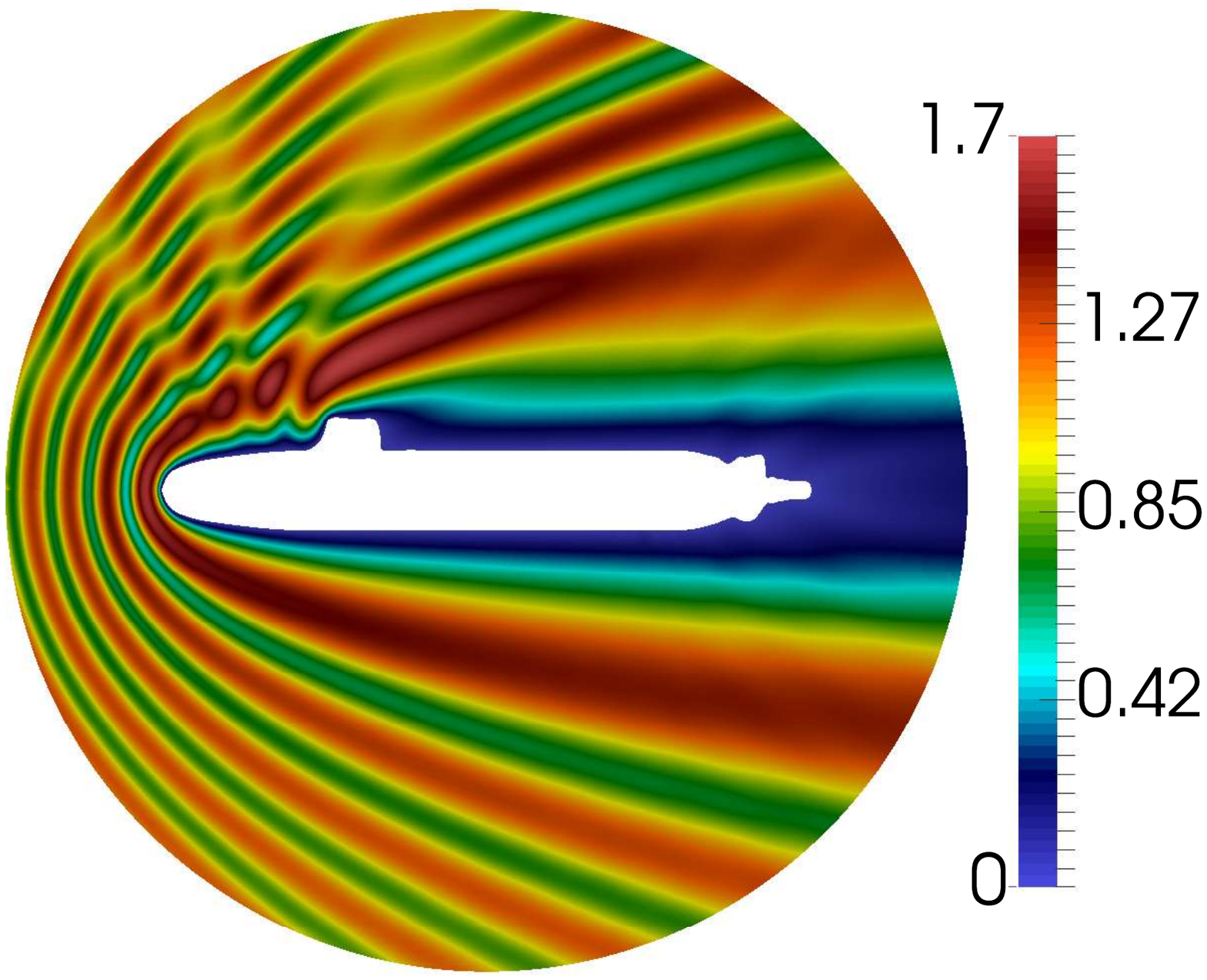} 
\caption{Boundary conforming grid for the computational domain between the prototype submarine and the artificial boundary at $R=2$ (left), and the magnitude of the total field (right)  for $k=20$, $p=5$, $ n_{\lambda}=5$, and $NT=5$. } 
\label{submesh}
\end{center}
\end{figure}
\begin{figure}[!h]
\begin{center}
\includegraphics[width=0.48\textwidth]{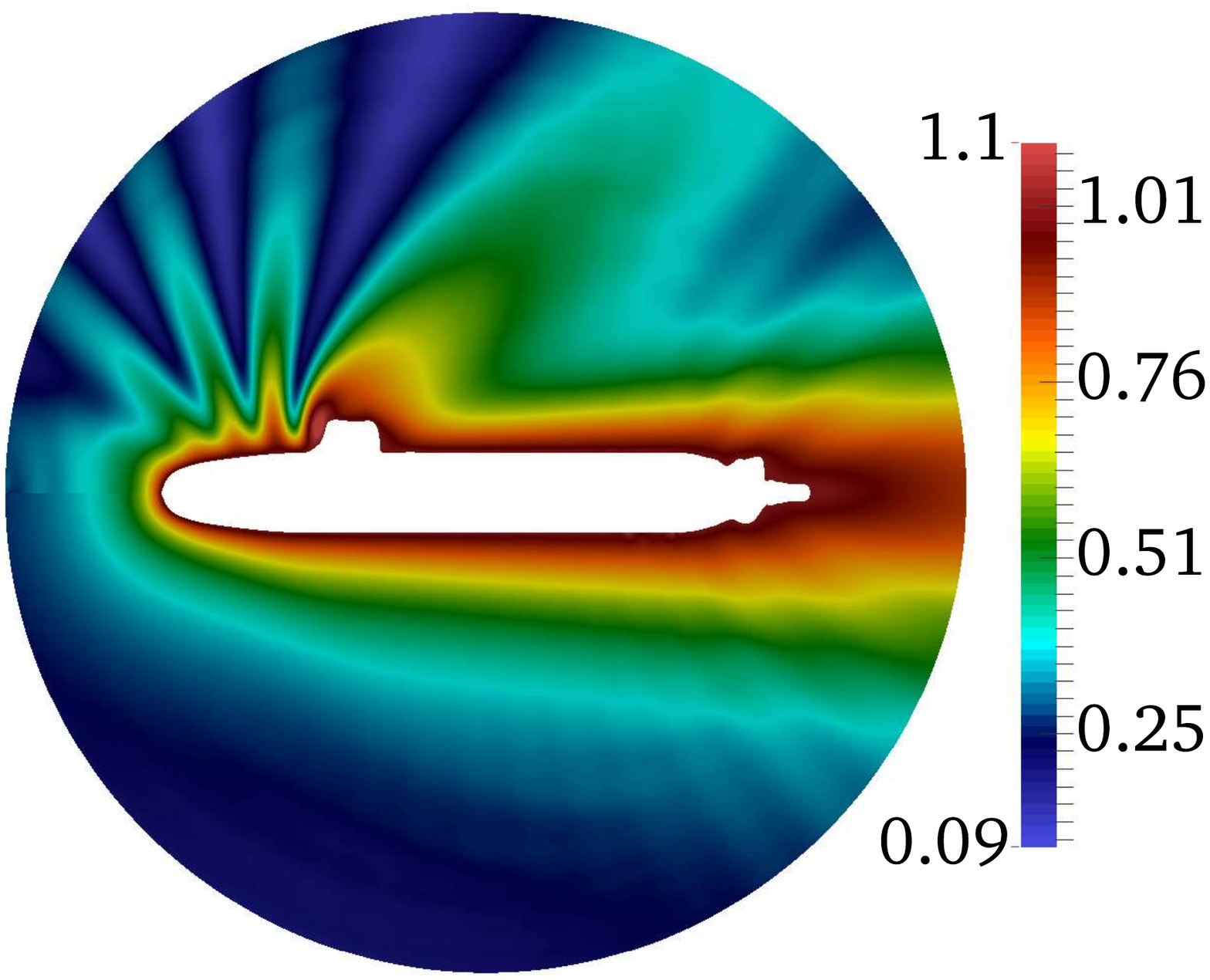} 
\includegraphics[width=0.5\textwidth]{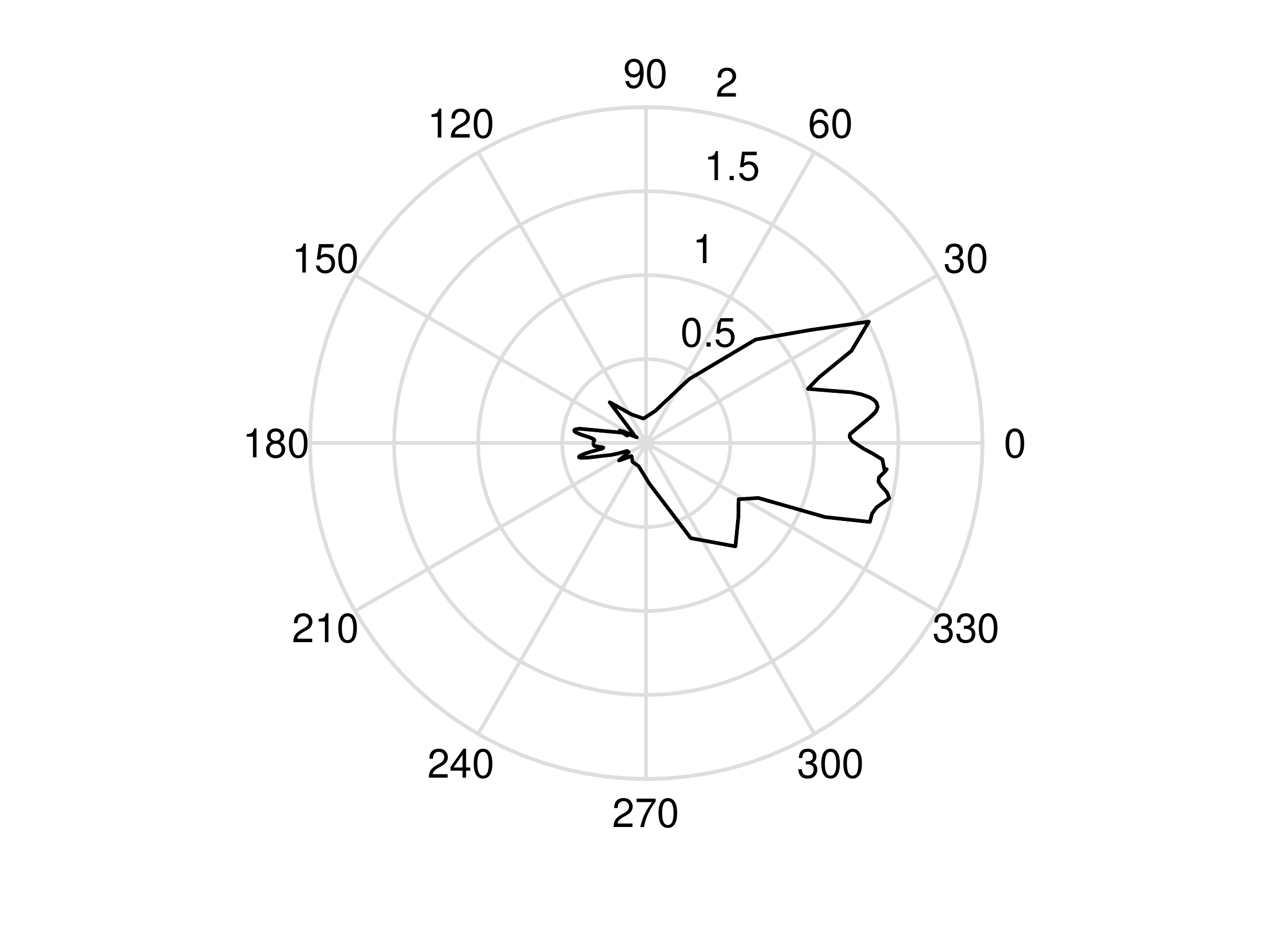}
\caption{The magnitude of the prototype submarine scattered field (left) and its corresponding farfield pattern (right) for $k=20$, $p=5$, $ n_{\lambda}=5$, and $NT=5$.} 
\label{Sub_k20}
\end{center}
\end{figure}
\begin{figure}[!h]
\begin{center}
\includegraphics[width=0.7\textwidth]{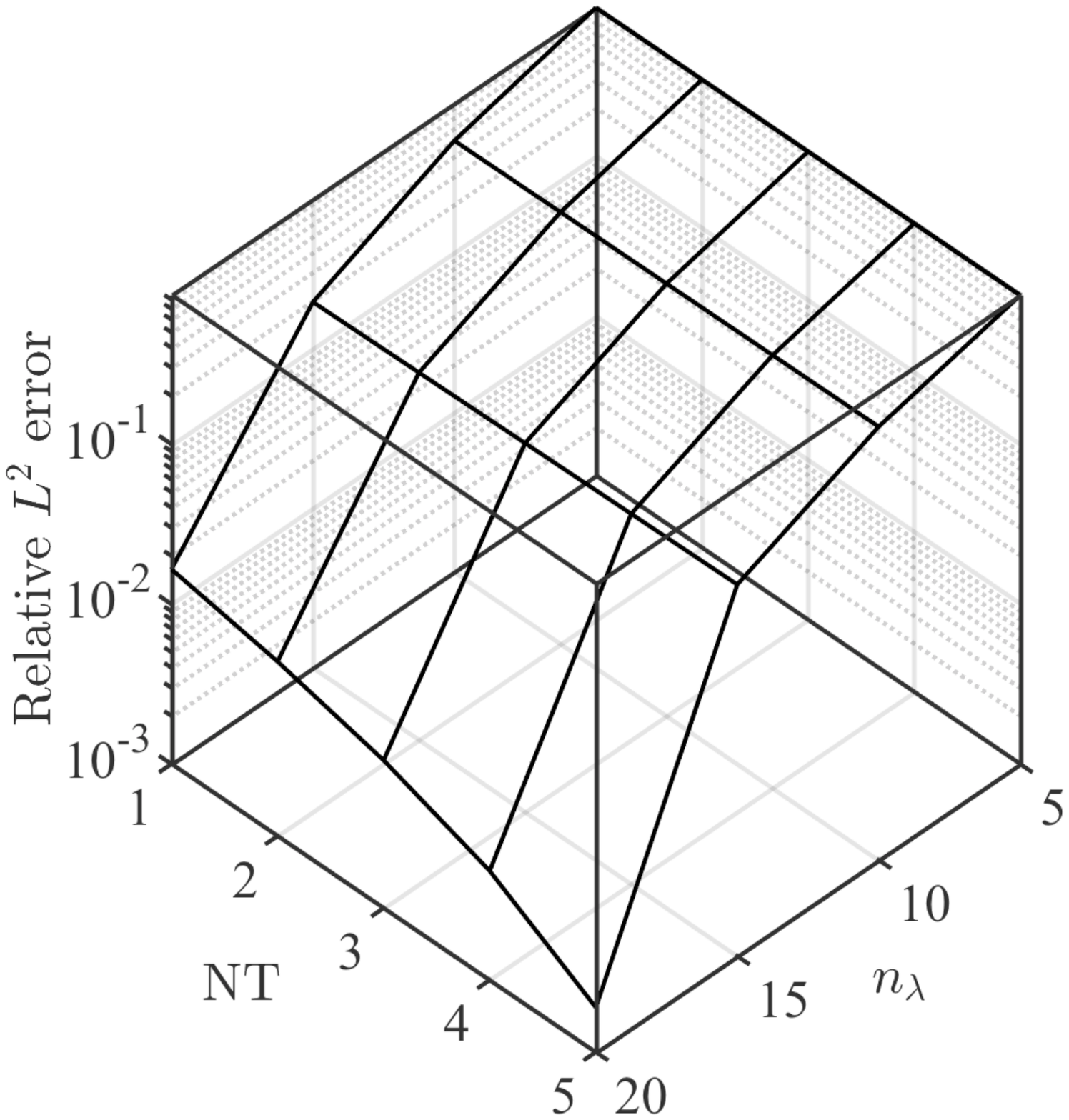}  
\caption{The FFP relative $L^2$ error for the prototype submarine employing a basis of order $p=5$ with wavenumber $k=20$. } 
\label{subconvk20}
\end{center}
\end{figure}

The Karp farfield expansion absorbing boundary condition is placed on an artificial boundary defined by a circle of radius $R=2$. The wavenumber, the basis order, the discretization density, and the Karp expansion number of terms are given by 
$k=20$, $p=5$, $ n_{\lambda}=5$, and $NT=5$, respectively.
The submarine mesh is obtained by adopting a high quality elliptic grid generation method \cite{ETNA09} for B-Spline mesh generation. As a consequence, the grid is smooth and it conforms well to the complex submarine boundary. This is illustrated in Fig. \ref{submesh} (left). 
Also, the numerical solution of the total field is shown to the right of Fig. \ref{submesh} while the scattered field and its corresponding farfield pattern are shown in Fig. \ref{Sub_k20}. As expected, the highest magnitude for the total field is located at the submarine conning tower. 

Since there is no exact solution available for this example, we calculate the relative $L^2$ error made in the computation of the FFP by considering an overly refined numerical solution as the reference solution. 
The dependence of this relative $L^2$ error  on $n_{\lambda}$, and $NT$ is illustrated in the surface graph depicted in Fig. \ref{subconvk20} for wavenumber $k = 20$, and basis order $p=5$. We observe that it decreases as both the discretization density $n_{\lambda}$ and the number of terms $NT$ of the Karp's expansion  increase. This behavior is completely analogous to the one observed for the circular cylinder example in the previous sections. 

\subsection{Acoustic scattering from a circular cylinder at high frequencies}
\label{HighFreq}
\begin{figure}[!h]
\begin{center}
\includegraphics[width=0.6\linewidth]{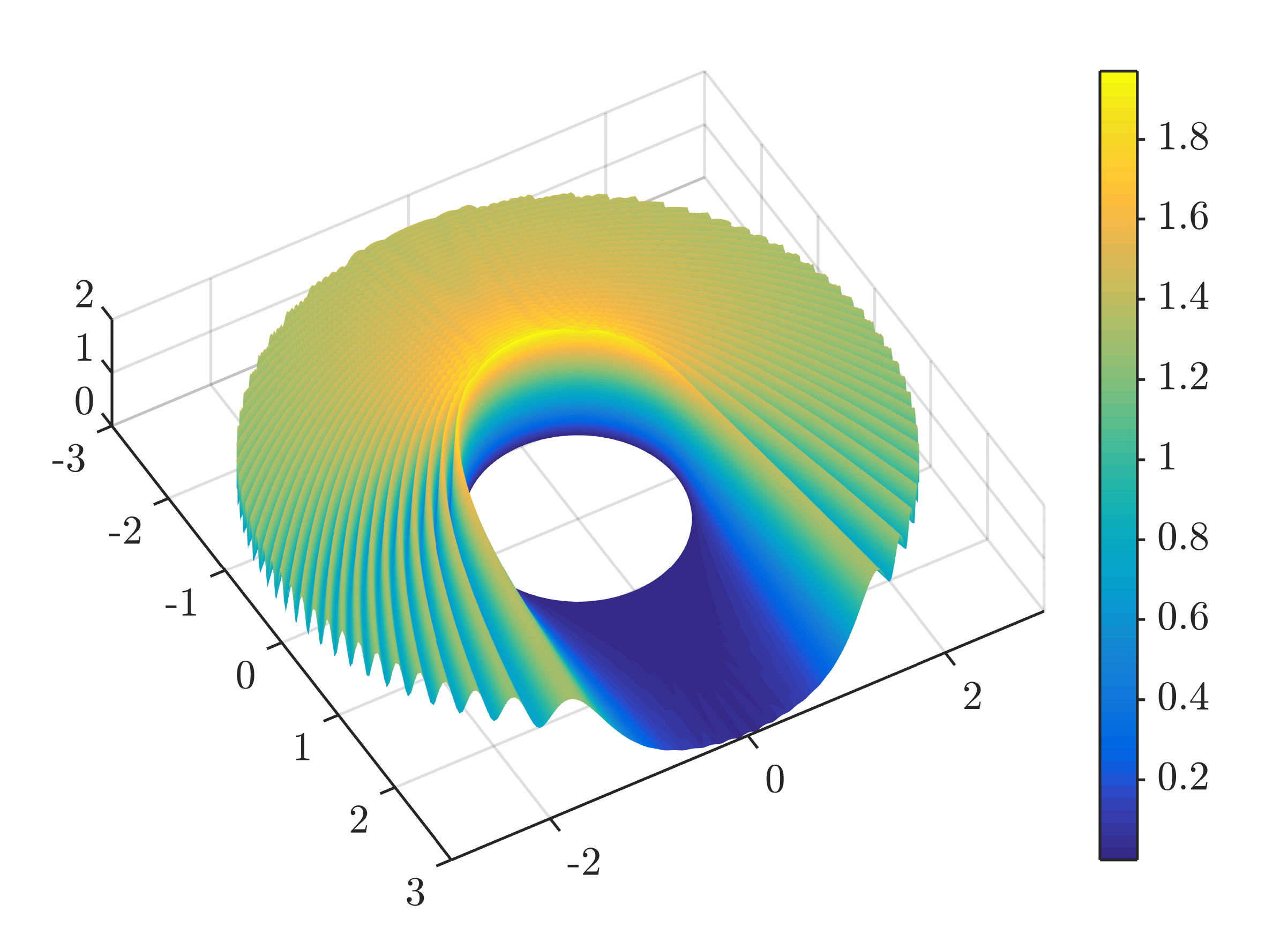}  
\caption{Total acoustic field for $R=3$, $p=5$, $k=50$, and $NT=5$.} 
\label{TotalFieldk_50}
\end{center}
\end{figure}

One of the most important finding of the application of IGA-FEABC to acoustic scattering problems is its ability to approximate the scattered field at very high frequencies such as $k=350$. We performed a series of experiments for the sound-soft circular cylindrical satterer increasing the frequency $k$. First, we analyzed the case for $k=50$. The total scattered field is illustrated in Fig. \ref{TotalFieldk_50} for $R=3$, $p=5$, and $NT=5$. 
\begin{figure}[!ht]
\begin{center}
\includegraphics[width=0.45\textwidth]{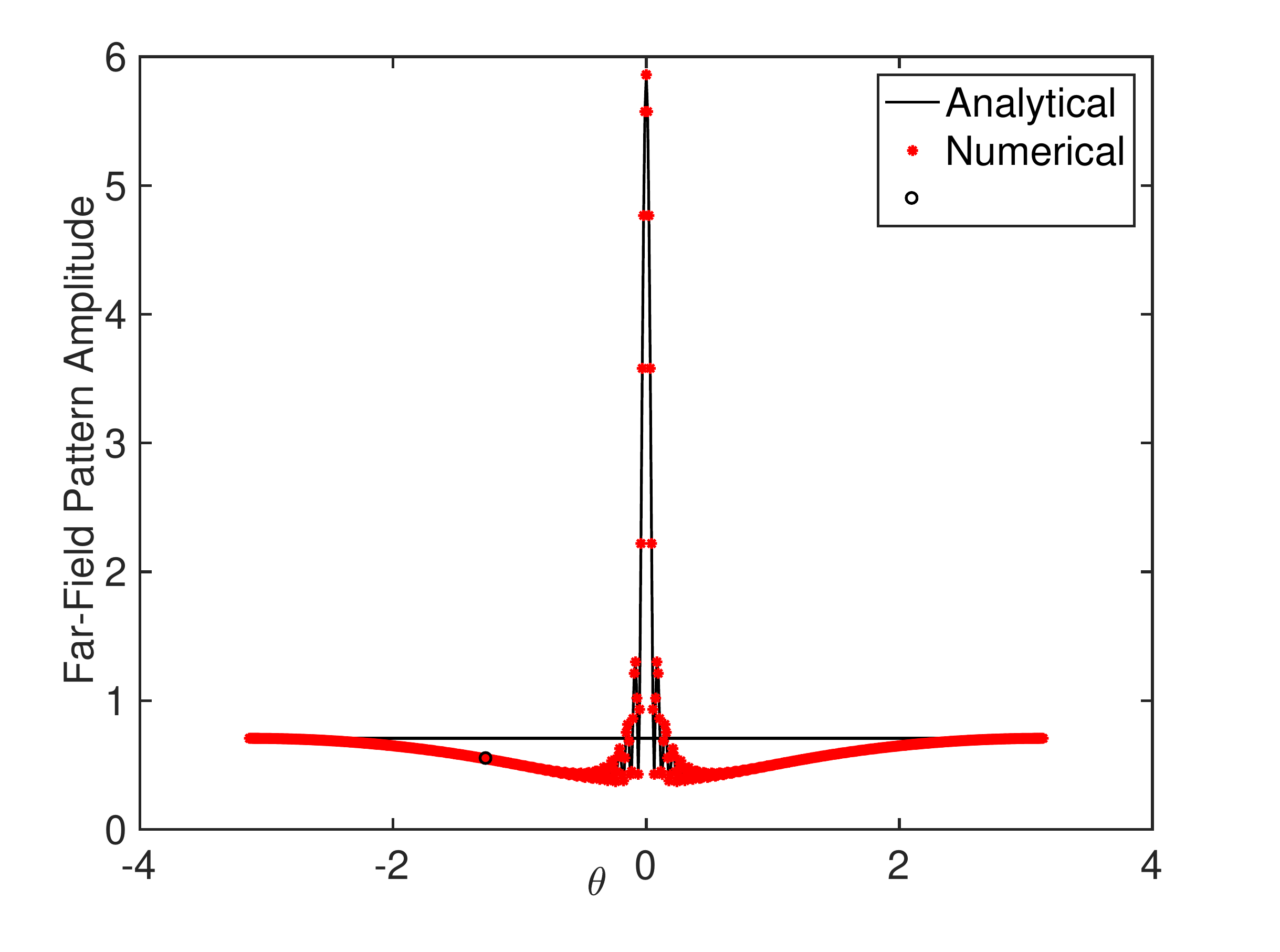}
\includegraphics[width=0.44\textwidth]{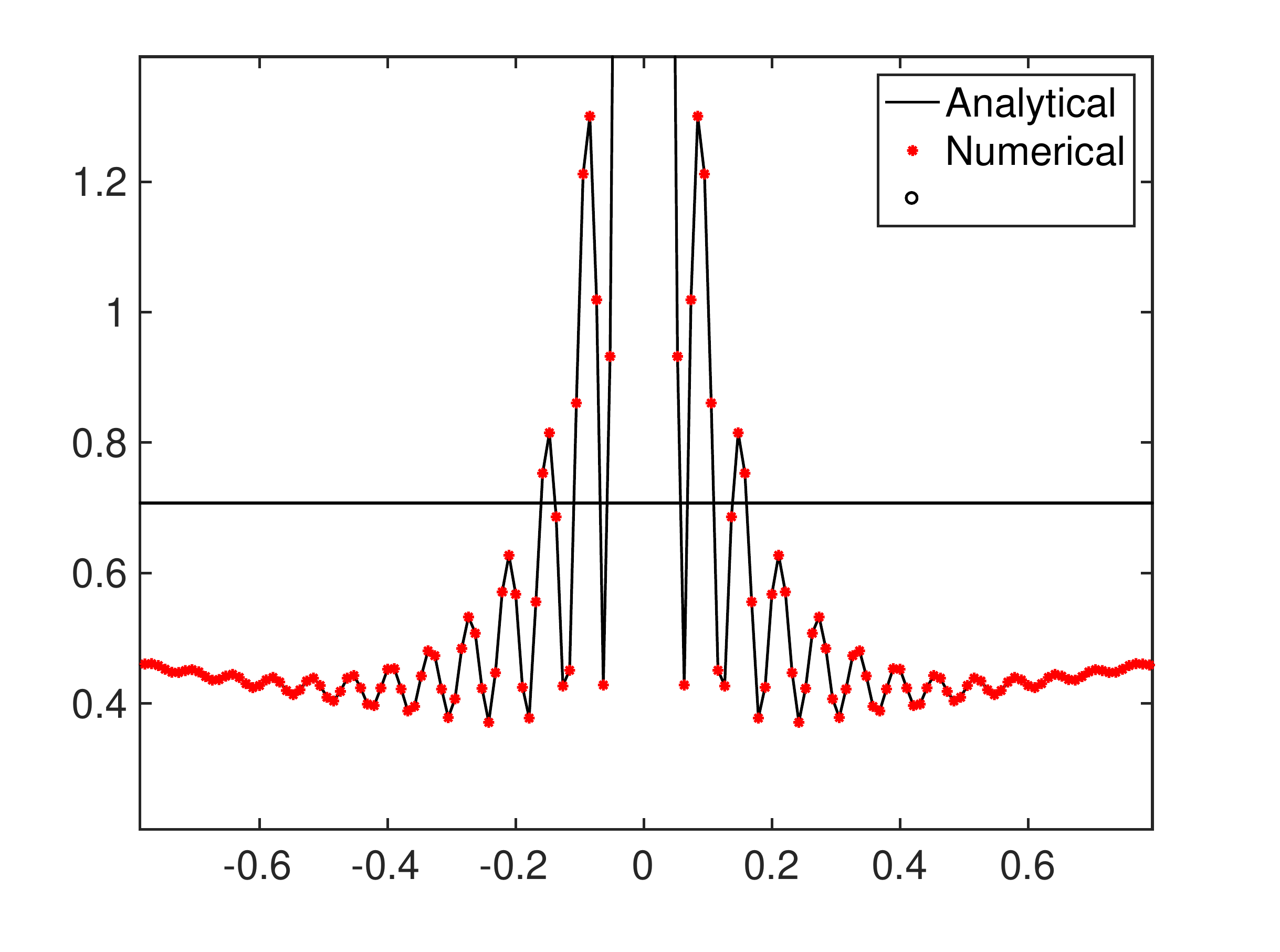}
\caption{Comparison of the numerical FFP versus the exact FFP for $k=50$, $p=5$, $NT=5$,  and $R=3$.}
\label{FFPK50}
\end{center}
\end{figure}
\begin{figure}[!ht]
\begin{center}
\includegraphics[width=0.45\textwidth]{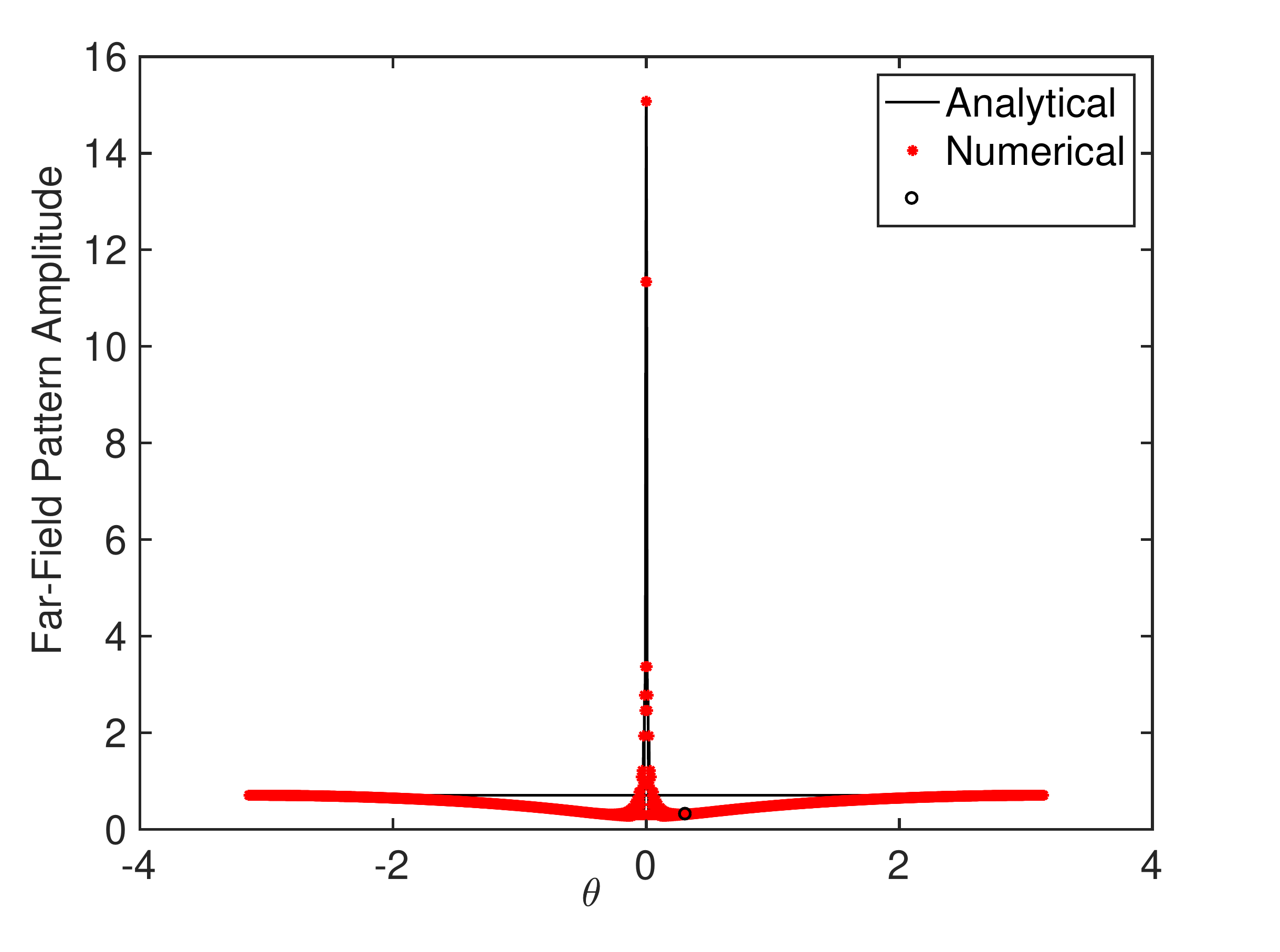}
\includegraphics[width=0.45\textwidth]{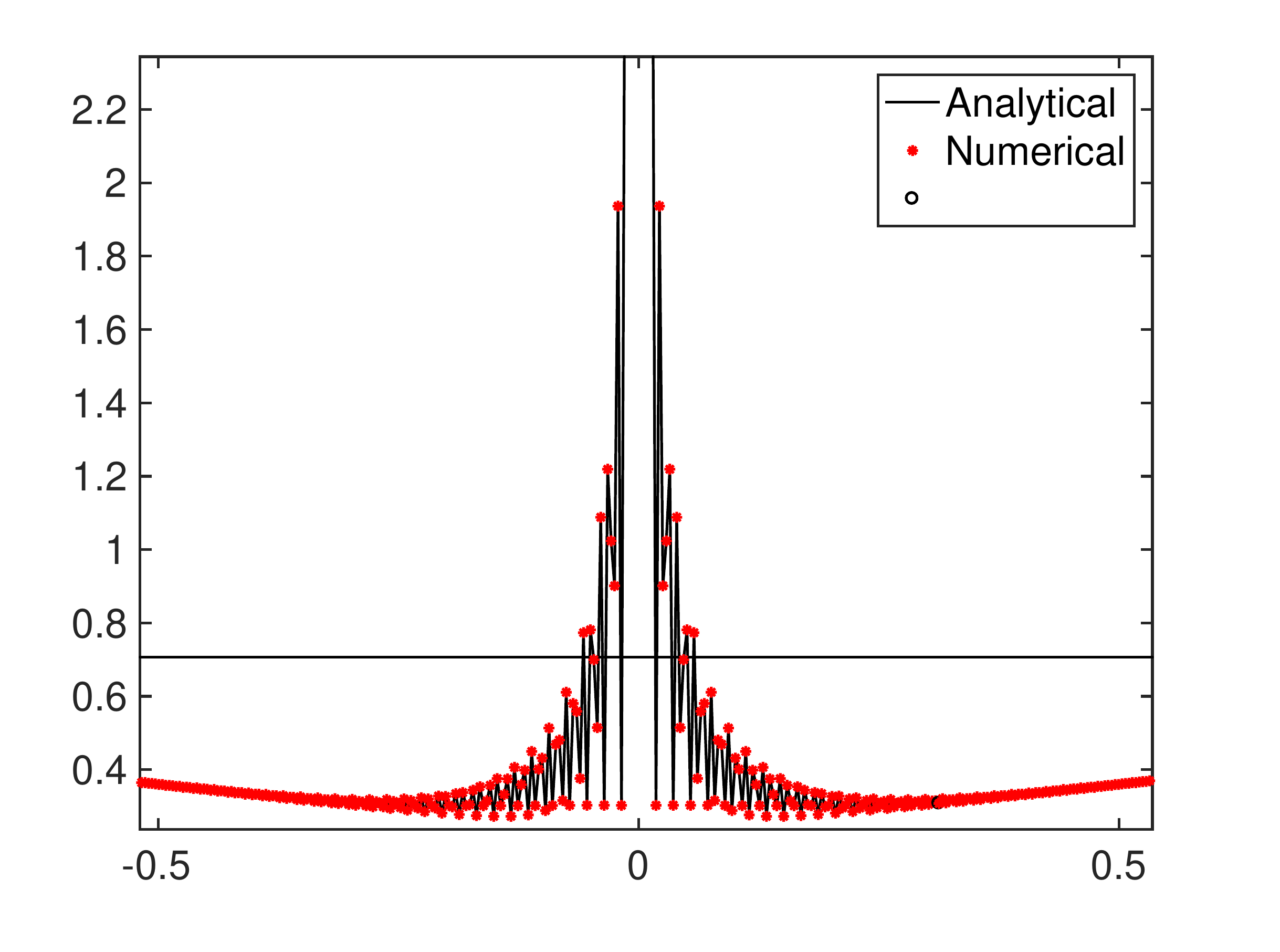}
\caption{Comparison of the numerical FFP versus the exact FFP for $k=350$, $p=5$, $NT=1$,  and $R=3$.}
\label{FFPK350}
\end{center}
\end{figure}

It is worth noticing how the IGA-FEABC capture the rapid oscillations in the high frequency regime. This is clearly shown in Fig. \ref{FFPK50} for $k=50$ and in Fig. \ref{FFPK350} for $k=350$. These figures illustrate the farfield patterns for these two high frequencies. Due to the high wave number, it is not possible to visualize the oscillations on the full domain $[-\pi,\pi]$ exhibited in the left windows of Figs. \ref{FFPK50} and \ref{FFPK350}. However by zooming these figures in the neighborhood of $\theta=0$ (right windows), it is possible to observe 
how the numerical solutions accurately adjust to the rapid oscillations.  The relative $L^2$ error of the farfield pattern approximation is close to $10^{-4}$ in both cases.
The remarkable fact is that these small errors are obtained with discretization densities $\nl=12$ for $k=50$ and $\nl=5$ for $k=350$.  

\begin{figure}[!h]
\begin{center}
\includegraphics[width=0.45\textwidth]{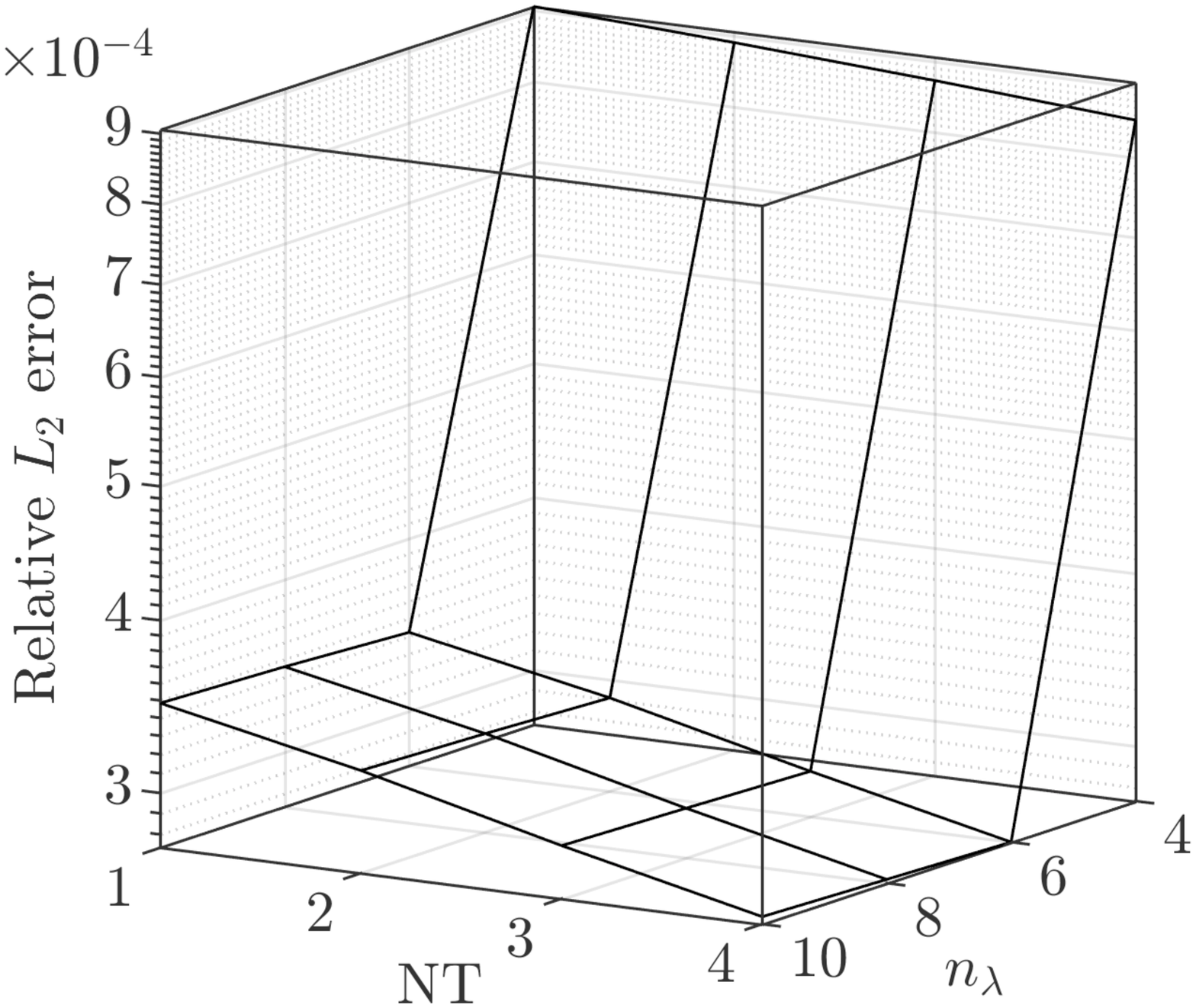} 
\includegraphics[width=0.45\textwidth]{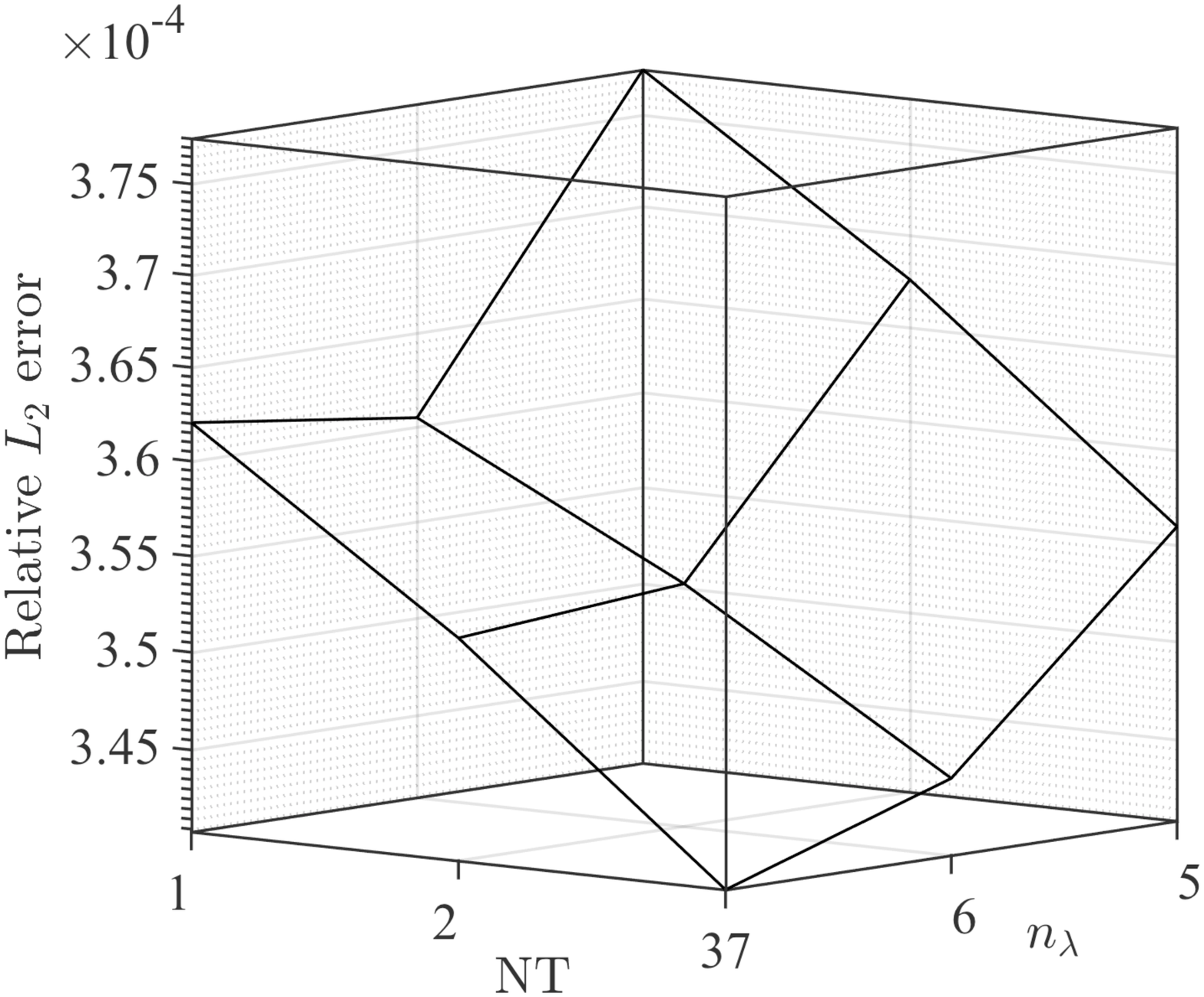} 
\caption{Farfield Pattern relative $L^2$ error for $p=5$ in terms of $NT$ and $\nl$ for $k=50$ (left) and
$k=100$ (right).} 
\label{SurfHF}
\end{center}
\end{figure}

 It is well known that conventional FEM may suffer from accuracy degradation in higher modes and Gibbs phenomena in wave propagation due to optical branches. A study of Isogeometric Analysis showed that it is possible to eliminate the optical branches of frequency spectra through nonlinear parameterization of the geometrical mapping 
 \cite{Hughes2014}.
 It is also known that  IGA NURBS bases outperforms conventional FEM in wave propagation analysis providing higher accuracy per degree of freedom and less dispersion error \cite{Hughes2008,Cotrell2006}.
The surface graphs in Fig. \ref{SurfHF} for $k=50$ and $k=100$ show the dependence of the $L^2$ relative error for $p=5$ fixed and $R=3$ with respect to $NT$ and $\nl$.  
We note that the behavior for high frequencies is similar to those observed for moderate frequencies. In fact, the relevant parameters responsible to reduce the error are the number of terms $NT$ of the KFE, the discretization density $\nl$, and the order $p$ of the basis. The minimum relative $L^2$ error obtained is close to $10^{-4}$ in both cases. The most remarkable fact of our results is that errors that low are still maintained for frequencies as high as $k=350$, as shown in Fig. \ref{FFPK350}. This is due to negligible pollution error in IGA of order $p=3$ and higher as observed in  \cite{Tahsin2017,Tahsin2016}.

\begin{figure}[!h]
\begin{center}
\includegraphics[width=0.6\textwidth]{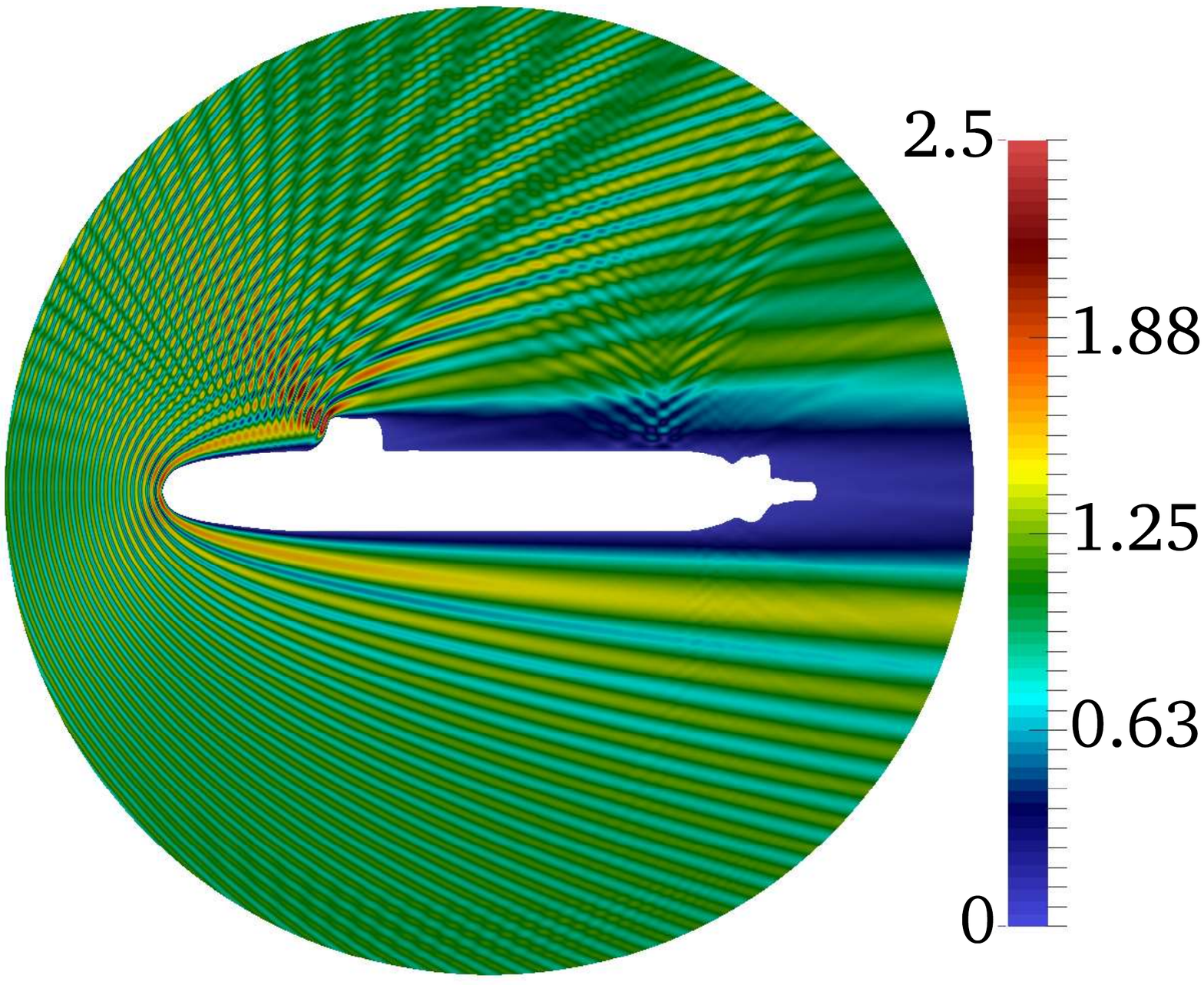}
\caption{Amplitude of the total field for the prototype submarine with $k=100$, $p=5$, $ n_{\lambda}=10$, and $NT=5$.} 
\label{Sub_k100}
\end{center}
\end{figure}
\begin{figure}[!h]
\begin{center}
\includegraphics[width=0.5\textwidth]{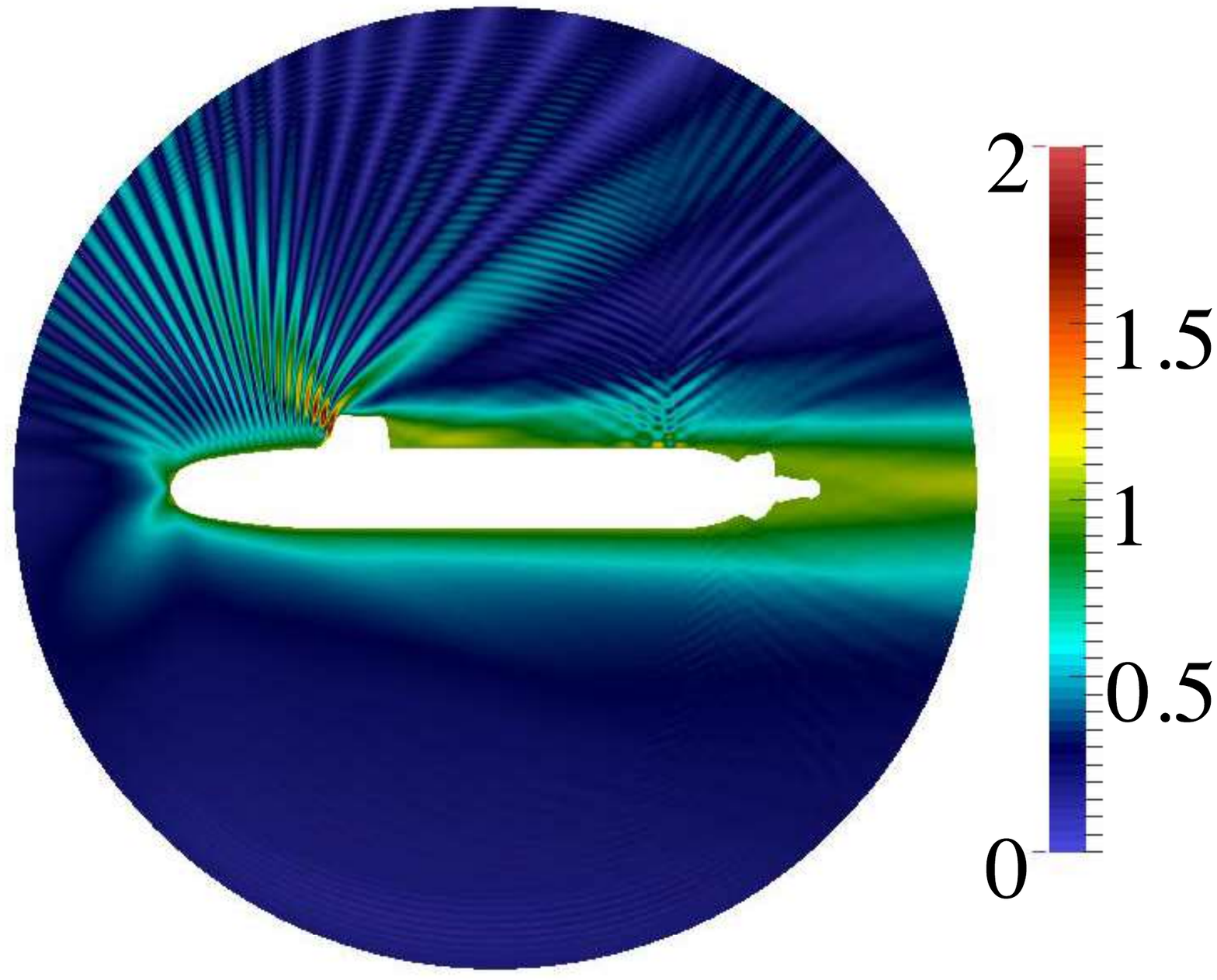} 
\includegraphics[width=.42\textwidth]{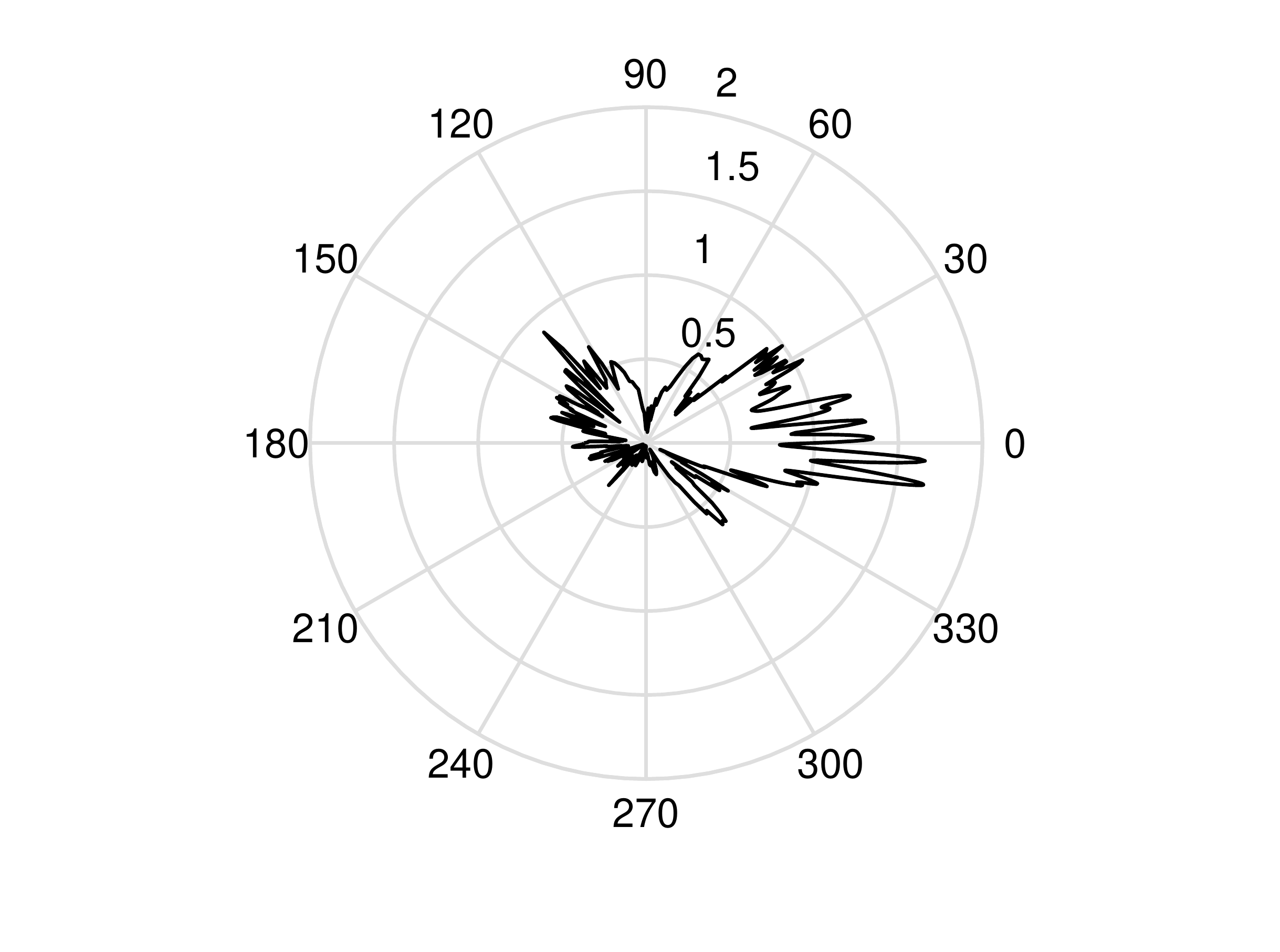}  
\caption{Amplitude of scattered field for the prototype submarine (left) and its farfield pattern for $k=100$, $p=5$, $n_{\lambda}=5$, $NT=5$.}
\label{SubFFP100}
\end{center}
\end{figure}

We also perform experiments for the prototype submarine of the previous section with wavenumber $k=100$. The magnitude of the scattered field and the total field are shown in Fig.\ref{Sub_k100} where $p=5$, $n_{\lambda}=10$, $NT=5$. The corresponding farfield pattern is also shown in Fig.\ref{SubFFP100}. The IGA-FEABC method produces numerical results qualitatively correct, similar to those obtained for $k=20$. The anticipated difference consists of the increment of the oscillations which are well-captured by the proposed numerical method.

\subsection{Plane wave scattering from a sphere.}
\label{Sphere}
\begin{figure}[!h]
\begin{center}
\includegraphics[width=0.48\textwidth]{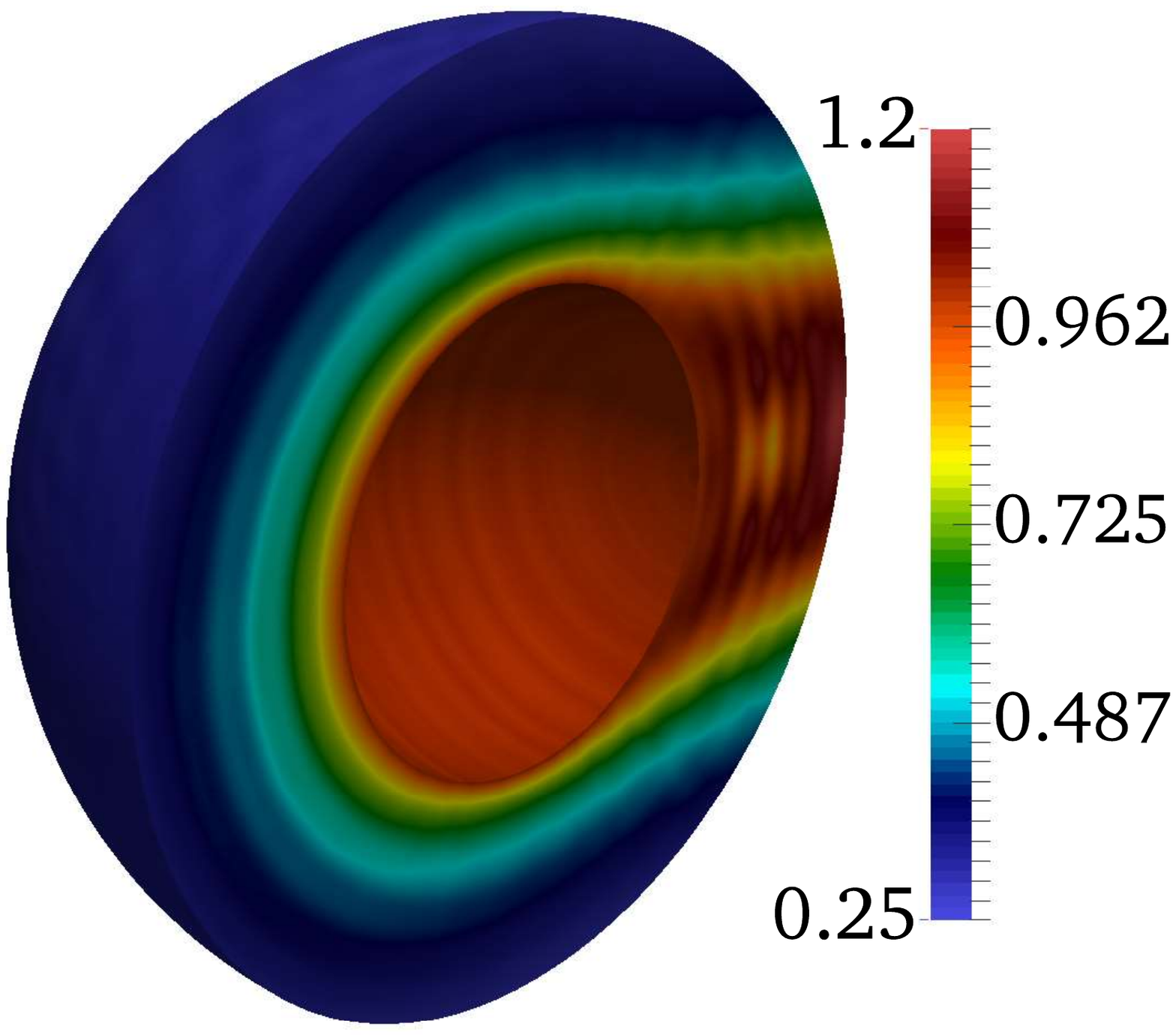} 
\includegraphics[width=0.48\textwidth]{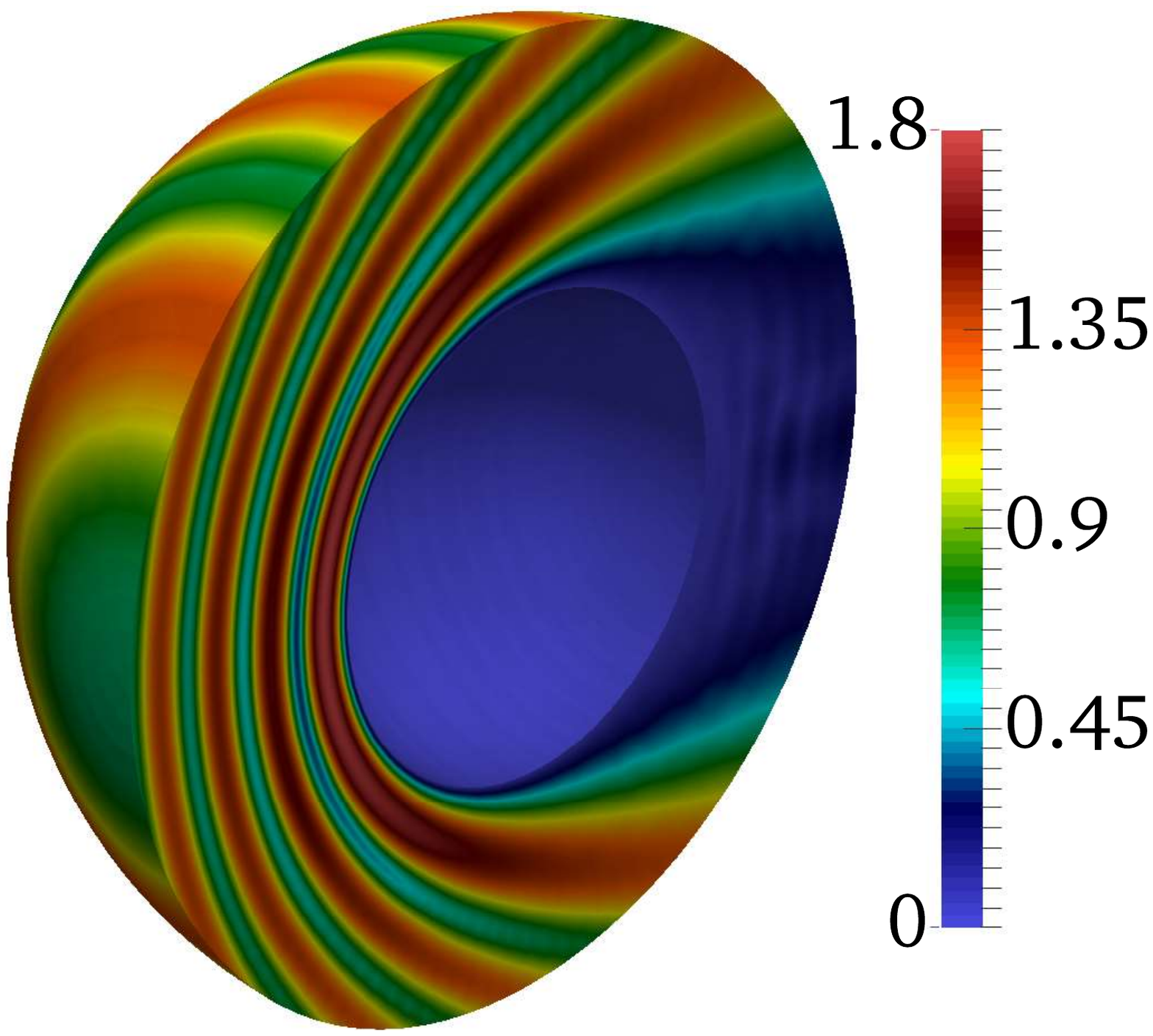} 
\caption{Cross-section of the sphere showing the scattered field amplitude (left) and total field amplitude (right) for $k=4\pi$, $p=5$, $\nl = 5$, and $NT=10$. }
\label{Sphere2pi}
\end{center}
\end{figure}

Finally, we analyze the numerical approximation for the scattering of a plane wave propagating in the positive direction of the $z$-axis from a spherical scatterer. The mathematical model  in weak form given by equations (\ref{WF3D1})-(\ref{WF3D3}) was formulated in Section \ref{Weakformulation3D} . 
This problem is axisymmetric about the $z-$axis. Therefore, the governing equations for the approximation $u$ of the scattered field is independent of the polar angle $\phi$. 
The angular coefficients $F_l$ of the Wilcox farfield expansion (WFE) are also independent of $\phi$. The amplitude of the scattered, and total fields  for a plane wave scattering from a sound soft sphere of radius $R_s=1$ are illustrated in 
 Fig. \ref{Sphere2pi}. They are  depicted 
 on a cross section $\phi$=constant of the computational domain. In this experiment, the artificial boundary consists of a sphere of radius $R=2$, the wavenumber is  $k=4 \pi$, the order of the NURBS basis used is $p=5$, the discretization density $n_{\lambda}=5$, and the number of terms in the Wilcox's expansion is $NT=10$.  
\begin{figure}[!h]
\begin{center}
\hspace{-1.7cm}
\includegraphics[width=0.5\textwidth]{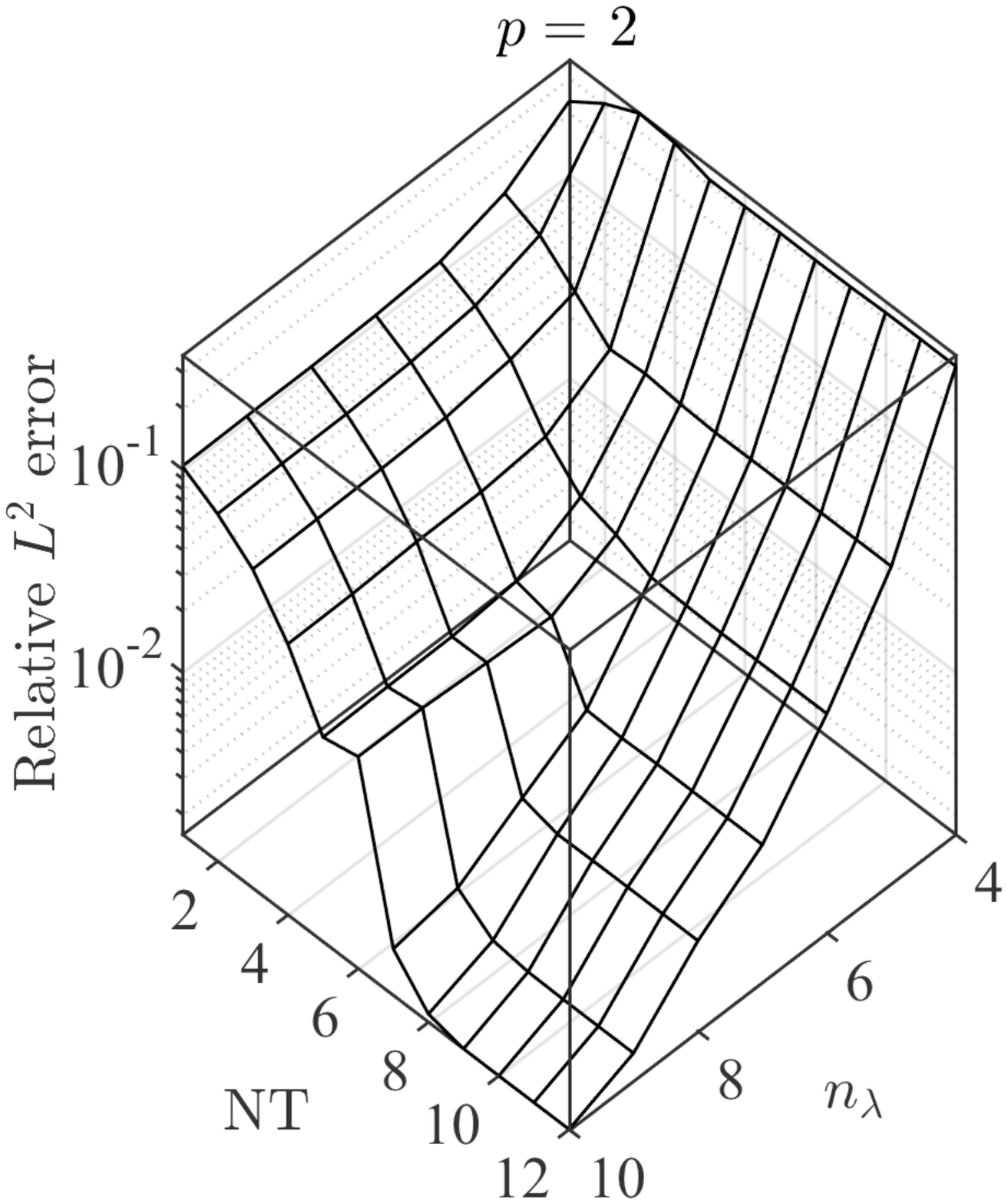} 
\vspace{-1.2cm}
\includegraphics[width=0.52\textwidth]{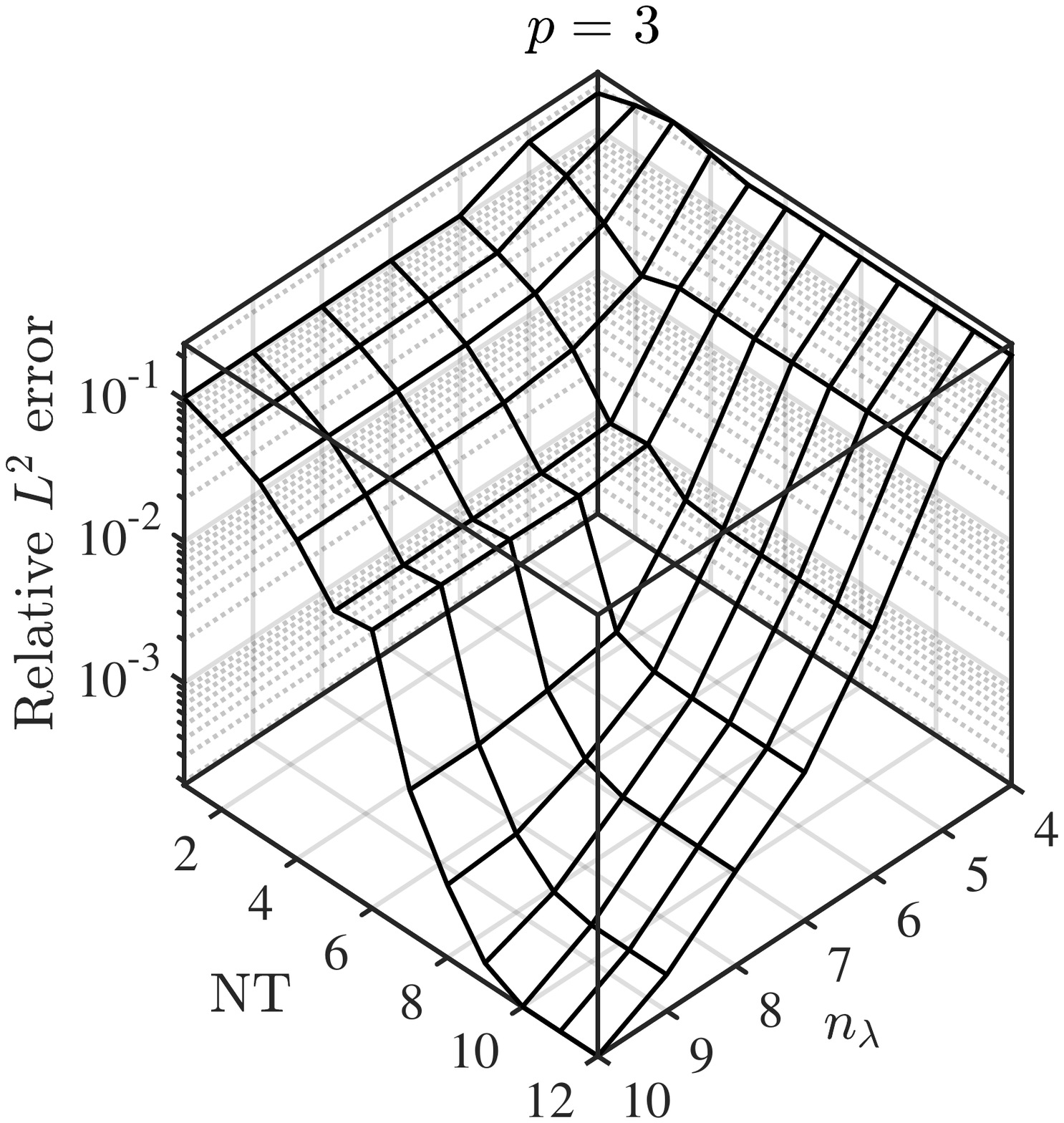} 
\includegraphics[width=0.55\textwidth]{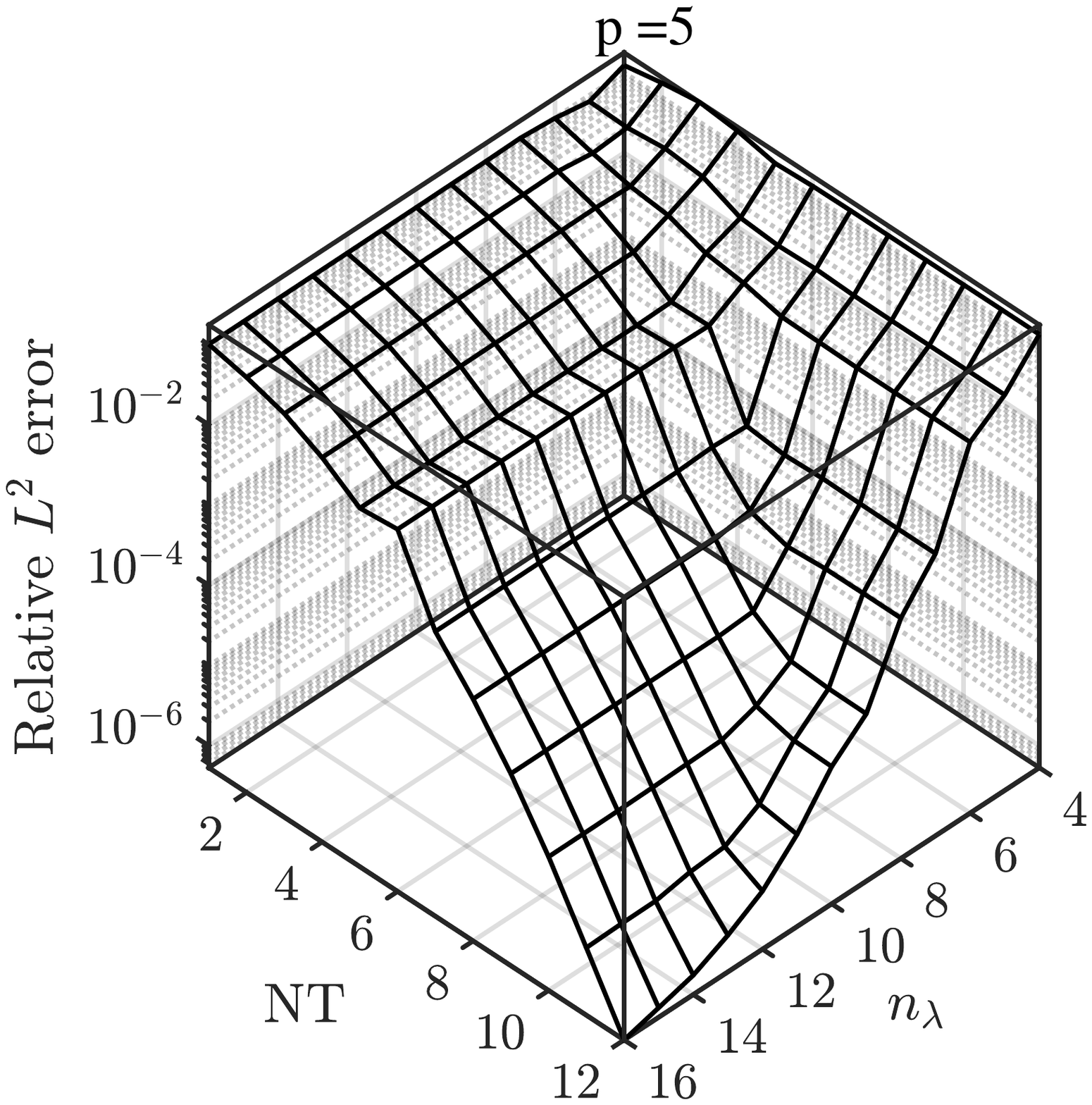} 
\caption{Relative $L^2$ error  at the artificial boundary for $k=2\pi$ under $h$-refinement varying $NT$ for various $p$.}
\label{Err3DK2pi}
\end{center}
\end{figure}

The surface graphs of Fig. \ref{Err3DK2pi} correspond to the relative $L^2$ error at the artificial boundary obtained for the same sphere scattering problem described above but  with frequency $k=2\pi$ instead.  
The analytical solution of this BVP can be obtained by eigenfunction expansions \cite{ColtonKress02}. We compare it against our numerical solution for the scattered field along the artificial boundary.

The 
relative $L^2$  error obtained by simultaneously performing $h$-refinement and increasing the number of terms $NT$ of the WFE is shown in Fig. \ref{Err3DK2pi}.
As previously observed in the 
 the 2D experiments, the $h$-refinement by itself is limited by the number of terms NT to decrease the error, similarly increasing $NT$ alone for a fixed discretization density $\nl$ is not enough to reduce the error. From the sequence of surface graphs for $p=2$, $p=3$, and $p=5$ of Fig. 
 \ref{Err3DK2pi}, we can observe how the error is reduced by adopting higher order basis functions. Actually,  
  the combined effect of $p=5$, $NT=12$ with $\nl=16$ leads to a minimum error approximately equal to $7\times 10^{-7}$. 

\section{Conclusion}

We have developed a numerical method that coupled a local high order absorbing boundary condition (FEABC) with isogeometric analysis for acoustic scattering problems. The FEABC is defined from the Karp (in 2D) and Wilcox (in 3D) expansions. These expansions which are exact representation of the outgoing waves outside the artificial boundaries were used to bound the infinite physical domains. They need to be truncated for computational purposes. As a consequence, they represent the scattered wave exactly up to the truncation number $NT$. On the other hand, the isogeometric analysis technique
unites the powers of finite element methods to solve partial differential equations with the  accuracy of computer aided design (CAD)  in representing complex shapes \cite{Hughes2005,Hughes2009}.
 
The order of convergence and the accuracy of the approximation is controlled by the order $p$ of the  basis employed by the IGA, the discretization density $n_{\lambda}$ of the control points, and the number of terms $NT$ of the FEABC. Highly accurate results with relative errors at the level of machine precision, in some of our experiments, can be obtained by implementing $p$- and $h$-refinement, and using an appropriate number $NT$ of terms for the FEABC. Our numerical experiments included plane wave scattering from an infinite circular cylinder, acoustic scattering from a prototype submarine in 2D, and scattering from a spherical scatterer.
 
We want to highlight the results of three of our experiments. Firstly, the fidelity of the numerical farfield pattern to follow the rapid oscillations of the exact FFP in the high frequency regime. This was observed even for a frequency $k=350$ (see Fig. \ref{FFPK350} ) where the $L^2$ relative error was about $10^{-4}$ using only a discretization density $\nl=5$, and one term, $NT=1$, of KFE. Secondly, the highly accurate approximation achieved for very low frequencies such as $k=0.01$. In fact, as it is shown in Table \ref{LowFreq}, relative errors close to machine precision were obtained with a rather coarse mesh of $10\times 60$ elements and employing only $NT=3$ terms for the KFE absorbing boundary condition. Finally, we want to point out the extraordinary accuracy reached by the IGA-KFE method when the artificial boundary is extremely close to the scatterer boundary (see Fig. \ref{R105}). In fact, we were able to obtain a relative error $9.99669\times 10^{-16}$ for a circular scatterer of radius $r_0=1$ when the radius of the circular artificial boundary was $R=1.05$. The ultimate linear system corresponding to this experiment was solved using the default  direct solver of MATLAB R2017a running on the same computer specified in Section \ref{CircularClose}. The CPU time employed for this computation was only 11.67 seconds. This remarkable performance of the proposed method in such a small computational domain verifies the exceptional efficiency of the IGA-FEABC for acoustic scattering.
We are currently working or plan to work on the extensions of IGA-FEABC to acoustic multiple scattering, to time-domain scattering (wave equation), and to elastic waves.

\appendixtitleon
\appendixtitletocon
\begin{appendices}
\label{Appendix}
\section{ }
\noindent{\it A brief overview of Isogeometric Analysis}
\medskip

In IGA, both the physical model and the solution space are constructed by {\emph{B-spline/NURBS}} functions. B-Splines are defined as parametric functions mapping the underlying parametric space $\hat{\Omega} \subset \mathbb{R} ^{d_p}$ to physical space $\Omega \subset \mathbb{R}^d$. The B-Spline curves and surfaces are constructed as a tensor product of a set of control points and {\emph{ knot vectors}} in each spatial direction. A knot vector is a set of non-decreasing parametric coordinates and represented by $ \vect{\xi}=\{\xi_1,\xi_2,...,\xi_{(n+p+1)}\}, \xi_i \leq \xi_{i+1} $, where $ \xi_i $ is the $ i^{th} $ knot and $ i$ is the knot index, $ i=1,2,\cdots,n+p+1$, where $p$ is the polynomial order and $ n $ is the number of basis functions. An element in IGA can be defined as the mapping of a non-zero knot span from parametric to physical space. A knot $\xi_i$ has {\it multiplicity} $k_i$ when it is repeated $k$ times in the knot vector where the B-Spline exhibit $C^{p-k}$ continuity in the corresponding physical point.  Usually {\it open knot vectors} are used in IGA where the first and the last knot values appear $ p+1 $ times, in other words, first and last knots have $k=p+1$ multiplicity which results in $C^0$ continuity at the patch boundaries. Hence, it is possible to conveniently achieve high order continuity and therefore increase computational accuracy within the domain boundaries while satisfying the required $C^0$ continuity to solve Helmholtz equation everywhere. 
\begin{equation}
\boldsymbol{\xi} =\{ 
\underbrace{0, \cdots, 0}_\text{p+1 times}, \cdots, \underbrace{1, \cdots, 1}_{\text{p+1 times}} 
\}.
\end{equation}
B-Splines defined over open knot vectors are interpolatory at first and the last knots.  If the knots are spaced equally within the knot vector, the knot vector is called {\it uniform} otherwise it is {\it non-uniform}. 
A B-Spline basis function is defined recursively by the Cox-de Boor recursion formula starting with the zeroth order $ (p=0) $ basis function:
\begin{equation}
\hspace*{-6.5cm}  \text{for} \ p = 0, \ N_{i}^{p}(\xi) = 
\left \{ \begin{aligned}
    1 & \ \ \ \xi_i \leq \xi \leq \xi_{i+1},\\
    0 & \ \ \ \text{otherwise},                                                                                                                                                                                                                            \label{eq.Bsp1}
    \end{aligned}
\right.
\end{equation}
\begin{equation}
\hspace*{-.254cm} \text{for}\ p = 1, 2, 3, \cdots \ \ \ N_{i}^p(\xi) = \frac{ \xi - \xi_i}{\xi_{i+p}-\xi_i}N_{i}^{p-1}(\xi) + \frac{ \xi_{i+p+1} - \xi}{\xi_{i+p+1}-\xi_{i+1}}N_{i+1}^{p-1}(\xi),          \label{eq.Bsp2}
\end{equation}
where $0/0 $  is defined to be zero. The first-order B-Spline functions are identical to their Lagrangian (FEM) counterparts. B-Splines also provide the partition of unity property, $\sum_{i=0}^n N_{i,p}(\xi) = 1$. The number of required shape functions for order $p$ analysis is $p+1$ shape functions; $ N_{i,p}(\xi) \neq 0 $ only when $\xi \in [ \xi_i, \xi_{i+p+1} ] $. Larger support in IGA translates into more expensive matrix assembly. Regardless, the total number of functions that any particular shape function can share support with is $2p+1$ in both IGA and conventional FEM. As a result, the linear matrix bandwidth in IGA is similar to that of FEM. Hence, solving the linear system in IGA is as expensive as that of conventional FEM for same order and number of degrees of freedom and the stiffness matrix is similarly sparse and symmetric.  In contrast to FEM, the shape functions in IGA are non-negative. Exact representation of both polynomials and conic sections such as circles, spheres, and ellipsoids can be generated using Non Uniform Rational B-Splines (NURBS). 

\begin{equation}
R_{i}^p(\xi) = \frac{N_{i}^{p}(\xi)w_i}{W(\xi)} = \frac{N_{i}^p(\xi)w_i}{\sum_{i=1}^{n} N_{i}^{p}(\xi) w_i},                                                                                                                                    \label{eq:NURBS}
\end{equation}
where $ \{N_i^{p} \}_{i=1}^n$ is a set of B-Spline basis functions and $ \{w_i \}_{i=1}^n $ is a set of positive NURBS weights.  If the weights are all equal, NURBS basis functions will reduce to their B-Spline counterparts, $R_{i}^{p} = N_{i}^{p}$, and the corresponding curve becomes a non-rational polynomial again. Hence, B-Splines are a subset of NURBS. Multivariate NURBS bases functions are generated as the tensor product of univariate basis:
\begin{equation}
R_{i,j}^{p,q} (\xi, \eta) = \frac{ N_{i}^{p}(\xi) M_{j,q}(\eta) w_{i,j}}{\sum_{i=1}^{n} \sum_{j=1}^{m} N_{i}^{p} (\xi) M_{j,q}(\eta) w_{i,j}},                                                                             \label{eq:NURBSsurf}
\end{equation}
\begin{equation}
R_{i,j,k}^{p,q,r} (\xi, \eta, \zeta) = \frac{ N_{i}^{p}(\xi) M_{j}^{q}(\eta) L_{k}^{r}(\zeta) w_{i,j,k}}{\sum_{i=1}^{n} \sum_{j=1}^{m}\sum_{k=1}^{l} N_{i}^{p} (\xi) M_{j}^{q}(\eta) L_{k}^{r}(\zeta) w_{i,j,k}}.                                                                      \label{eq:NURBSvol}
\end{equation}
where $ N_{i}^{p}(\xi)$, $M_{j}^{q}(\eta)$, and $ L_{k}^{r}(\zeta)$ are B-Spline basis functions of order $p$, $q$, and $r$ respectively. NURBS curves, surfaces, and volumes are defined as a linear combination of these basis functions and the corresponding control points denoted with $\vect{B}$:
\begin{equation}
\begin{aligned}
&\vect{C}(\xi) = \sum_{i=1}^{n} R_i^p(\xi) \vect{B}_i, \\
&\vect{S}(\xi, \eta) = \sum_{i=1}^{n} \sum_{j=1}^{m} R_{i,j}^{p,q}(\xi, \eta) \vect{B}_{i,j},\\                                                                                                                                                                                                                                                                                                                                                                                                                       
&\vect{V}(\xi, \eta) = \sum_{i=1}^{n} \sum_{j=1}^{m} \sum_{k=1}^{l}R_{i,j,k}^{p,q,r}(\xi, \eta, \zeta) \vect{B}_{i,j,k}.                                                                                                                                                                                                                                                                                                                                                                                                                                                           \label{eq:NURBScurve}
\end{aligned}
\end{equation}


\end{appendices}

%


\section*{Acknowledgments}
The second author acknowledges the support provided by the Office of Research and Creative Activities (ORCA) of Brigham Young University.


\bibliographystyle{elsarticle-num}
\bibliography{AcoBib}

\end{document}